\documentclass[11pt]{article}
\usepackage{times}
\usepackage{amsmath,amsfonts,amssymb,latexsym,epsfig}
\usepackage{color,epsf,float}
\usepackage{graphicx,stix}
\usepackage[toc]{appendix}

\usepackage{alltt}
\usepackage{url}
\usepackage{a4}

\newcommand{\e}{\ensuremath{\mathrm{e}}}
\newcommand{\ad}{\ensuremath{\mathrm{ad}}}

\begin{document}
\title{Exponential perturbative expansions and coordinate transformations}

\author{Ana Arnal\thanks{Email: \texttt{ana.arnal@uji.es}}
 \and
 Fernando Casas\thanks{Email: \texttt{Fernando.Casas@uji.es} }\and Cristina Chiralt\thanks{Email: \texttt{chiralt@uji.es}}}

\date{}

\maketitle

\begin{abstract}

We propose a unified approach for different exponential perturbation techniques used in the treatment of time-dependent quantum mechanical
problems, namely the Magnus expansion, the Floquet--Magnus expansion for periodic systems, the quantum averaging technique and the Lie--Deprit
perturbative algorithms.  Even the standard perturbation theory fits in this framework. The approach is based on carrying out an appropriate change of coordinates (or picture) in each case, 
and can be formulated for any time-dependent
linear system of ordinary differential equations. All the procedures  (except the standard perturbation theory) lead to approximate solutions preserving 
by construction unitarity when applied to the time-dependent Schr\"odinger equation.

\vspace*{1cm}

\begin{center}
Institut de Matem\`atiques i Aplicacions de Castell\'o (IMAC) and  De\-par\-ta\-ment de
Ma\-te\-m\`a\-ti\-ques, Universitat Jaume I,
  E-12071 Cas\-te\-ll\'on, Spain.
\end{center}

\end{abstract}\bigskip

\section{Introduction}
Linear differential equations of the form
\begin{equation}   \label{ecu.1}
\frac{dU}{dt} = A(t) \, U, \qquad U(0) = I
\end{equation}   
are ubiquitous in many branches of physics, chemistry and mathematics. 
Here $U$ is a real or complex $d \times d$ matrix, and $A(t)$ is a sufficiently smooth matrix to ensure the
existence of solutions. Perhaps the most important example corresponds to the Schr\"odinger equation for the evolution operator in quantum systems with
a time-dependent Hamiltonian $H(t)$, in which case $A(t) = -i H(t)/\hbar$. Particular cases include spin dynamics in magnetic resonance (Nuclear Magnetic
Resonance -NMR, Electronic Paramagnetic Resonance -EPR, Dynamic Nuclear Polarization -DNP, etc) \cite{slichter90pom,ernst86pon,mehring01oom},
electron-atom collisions in atomic physics, pressure broadening of rotational spectra in molecular physics, control of chemical reactions with driving 
induced by laser beams, etc. When the time-dependence of the Hamiltonian is periodic, as occurs for instance in 
periodically driven quantum systems, atomic quantum gases in periodically driven optical lattices, etc. \cite{eckardt17aqg,ota19feo}, the Floquet theorem  \cite{hale80ode} 
relates $H(t)$ with a constant
Hamiltonian. More specifically, it implies that the evolution operator is factorized as $U(t) = P(t) \exp(t F)$, with $P(t)$ a periodic time-dependent matrix
and $F$ a constant matrix. This theorem has been widely used in problems of solid state physics, in NMR, in the quantum simulation of systems with time-independent
Hamiltonians by periodically driven quantum systems, etc. \cite{mananga16otf}. The Average Hamiltonian Theory is also closely related with this result, and the effective 
Hamiltonian is an important tool in the description of the system \cite{maricq82aoa,goldman14pdq}.

In general, equation (\ref{ecu.1}) cannot be solved in closed form, and so different approaches have been proposed along the years to get approximations,
{both analytic and numerical. Among the former, we can mention the 
	standard perturbation theory, the average Hamiltonian theory, and the Magnus expansion \cite{ota19feo,goldman14pdq,casas01fte}. Concerning the second approach, different numerical algorithms
	have been applied to get solutions on specific time intervals \cite{lubich08fqt,blanes16aci,hairer06gni}.} 

{In this work we will concentrate on different techniques providing analytic approximations to the solution of (\ref{ecu.1}) that also share some of its most
	salient qualitative features. In particular, if (\ref{ecu.1}) represents the Schr\"odinger equation, it is well known that $U(t)$ is unitary for all $t$, and this
	guarantees that its elements represent probabilities of transition between the different states of the system. It happens, however, that not every approximate
	scheme (either analytic or numerical) renders unitary matrices, and thus the physical description they provide is unreliable, especially for large integration times.}

{The Magnus expansion \cite{magnus54ote} presents the remarkable feature that it allows one to express the solution as the exponential of a series,
	$U(t) = \exp(\Omega(t))$, so that even when the series $\Omega(t) = \sum_{n \ge 1} \Omega_n(t)$ is truncated, the corresponding approximation 
	is still unitary when applied to the Schr\"odinger equation. More generally, if (\ref{ecu.1}) is defined in a Lie group $\mathcal{G}$, then it provides an approximate
	solution also belonging to $\mathcal{G}$. Moreover, it has also been used to construct efficient numerical integration algorithms also preserving this feature
	\cite{iserles00lgm,blanes16aci}.  }

When $A(t)$ depends periodically on time with period $T$, the Magnus expansion does not explicitly provide the structure of the solution ensured by the
Floquet theorem, i.e., $U(t) = P(t) \e^{t F}$. with $P(t)$ periodic and $F$ a constant matrix. 
Nevertheless, it can be conveniently generalized to cover also this situation, in such a way that the matrix $P(t)$
is expressed as the exponential of a series of periodic terms. 
The resulting approach (the so-called Floquet--Magnus expansion \cite{casas01fte}) has been used during
the last years in a variety of physical problems  \cite{mananga16otf,eckardt17aqg,haga19dot,ota19feo}.

Very often, the coefficient matrix in (\ref{ecu.1}) is of the form $A(t) = A_0 + \varepsilon A_1(t)$, where $A_0$ is constant, $A_1(t+T) = A_1(t)$ and
$\varepsilon > 0$ is a (small) parameter. In other words, one is dealing with a  time-dependent 
perturbation of a solvable problem defined by $A_0$. In that case, several
perturbative procedures exist to construct $U(t)$ as a power series in $\varepsilon$, either directly (by applying standard perturbation techniques in the
interaction picture defined by $A_0$) or taking into account the structure ensured by the Floquet theorem and constructing both matrices $P(t)$ and $F$
as power series \cite{maricq82aoa}. Of course, if $P(t) = \exp(\Lambda(t))$ and $\Lambda(t)$ is constructed also as a power series in $\varepsilon$, then
the qualitative properties of the solution are inherited by the approximations (in particular, they are unitary in quantum evolution problems) \cite{casas12utd}.

An alternative manner of viewing the Floquet theorem is to interpret the matrix $P(t)$, provided it is invertible, as a time-periodic transformation to a new set of 
coordinates where the new evolution equation has a constant coefficient matrix $F$, so that $\exp(t F)$ is the exact solution in the new coordinates. In
the language of Quantum Mechanics, this corresponds to a change of picture.
This interpretation leads to the important mathematical notion of a \textit{reducible system}: according to Lyapunov, the general system (\ref{ecu.1}) is 
called reducible if there exists a matrix $P(t)$ which together with $\det(P^{-1}(t))$ is bounded on $0 \le t < +\infty$ such that the system obtained from
(\ref{ecu.1}) by applying the linear transformation defined by $P(t)$ has constant coefficients \cite{nemytskii89qto}. In this sense, if $A(t)$ is a periodic matrix,
then (\ref{ecu.1}) is reducible by means of a periodic matrix. The situation is not so clear, however, when $A(t)$ is quasi-periodic or almost periodic: 
although several general 
results exist ensuring reducibility \cite{jorba97ero,jorba92otr,afzal18rpf}, there are also examples of irreducible systems \cite{palmer80otr}.

In this work we pursue and generalize this interpretation to show that all the above mentioned exponential perturbative treatments can be considered as particular 
	instances of a generic change of variables applied to the original differential equation. The idea of making a coordinate change to analyze a problem
	arises of course in many application settings, ranging from canonical (or symplectic) transformations in classical mechanics (based either on 
	generating functions expressed in terms of mixed, old and new, variables \cite{arnold06mao}, or on the Lie-algebraic setting \cite{hori66tog,deprit69ctd}) to changes of picture and unitary transformations in quantum mechanics. What we intend here
	is to show that several widely used perturbative expansions in quantum mechanics can be indeed derived from the same basic principle using different
	variations of a unique ansatz based on a generic linear transformations of coordinates. We believe this interpretation sheds new light into the different
	expansions, and moreover allows one to elaborate a unique procedure for analyzing a given problem and compare in an easy way how they behave in practice.

It is important to remark
that all the procedures considered here (with the exception of course of the standard perturbation theory) preserve by construction the unitarity of the solution when (\ref{ecu.1}) refers to the Schr\"odinger equation. More
generally, the approximations obtained evolve in the same matrix Lie group as the exact solution of the differential equation (\ref{ecu.1}).

%%%%%%%%%%%%%%%%%%%%%%%%%%%%%%%%%%%%%%%%%%
\section{Coordinate transformations and linear systems}
\label{sec.2}

To begin with, let us consider the most general situation of a
linear differential equation
\begin{equation}   \label{ccv.1}
	\dot{x} \equiv \frac{dx}{dt} = A(t) x,  \qquad x(0) = x_0 \in \mathbb{C}^d,
\end{equation}
with $A(t)$ a $d \times d$ matrix whose entries are integrable functions of $t$. Notice that $U(t)$ in (\ref{ecu.1}) can be considered as the fundamental
matrix of (\ref{ccv.1}). We analyze 
a change of variables $x \longmapsto X$ transforming  the original system (\ref{ccv.1}) into
\begin{equation}  \label{ccv.2}
	\frac{dX}{dt} = F(t) X,  \qquad X(0) = X_0,
\end{equation}
where the matrix $F$ adopts some desirable form. Since we are interested in preserving qualitative properties of (\ref{ccv.1}), we impose
an additional requirement for the transformation, namely, it has to be of the form
\begin{equation}   \label{ccv.3}
	x(t) = \exp(\Omega(t)) X(t),  \qquad \mbox{ with } \qquad \Omega(0) = 0,
\end{equation}
so that $x(0) = X(0)$. Thus, in particular, if (\ref{ccv.1}) is the Schr\"odinger equation with Hamiltonian $H(t) = i \hbar A(t)$, $\Omega(t)$ is skew-Hermitian.

It is clear that the transformation (\ref{ccv.3}) is completely determined once the generator $\Omega(t)$ is obtained. An evolution equation for
$\Omega$ is obtained by introducing (\ref{ccv.3}) in (\ref{ccv.1}) and also taking into account (\ref{ccv.2}) as
\begin{equation} \label{eq.4a}
	\frac{d }{dt} \exp(\Omega) = A(t) \exp(\Omega) - \exp(\Omega) F(t).
\end{equation}   
The derivative of a matrix exponential can be written as \cite{blanes09tme}
\begin{equation}  \label{eq.4b}
	\frac{d }{dt} \exp(\Omega(t)) = d \exp_{\Omega(t)} (\dot{\Omega}(t)) \, \exp(\Omega),
\end{equation}
where the symbol $d \exp_{\Omega}(C)$ stands for the (everywhere convergent) power series
\begin{equation}  \label{eq.4c}
	d \exp_{\Omega}(C) = \sum_{k=0}^{\infty} \frac{1}{(k+1)!} \ad_{\Omega}^k (C) \equiv \frac{ \exp(\ad_{\Omega}) - I}{\ad_{\Omega}} (C).
\end{equation}
Here $\ad_{\Omega}^0 C = C$, $\ad_{\Omega}^k C = [\Omega, \ad_{\Omega}^{k-1} C]$,  and $[\Omega,C]$ denotes the usual commutator. 
By inserting (\ref{eq.4b}) into (\ref{eq.4a}) one gets
\begin{equation}   \label{ccv.4}
	d \exp_{\Omega} (\dot{\Omega}) = A - \exp(\Omega) F \exp(-\Omega) =   A -  \exp(\ad_{\Omega}) F,
\end{equation}
where
\[
\exp(\ad_{\Omega}) F = \sum_{n \ge 0} \frac{1}{n!} \ad_{\Omega}^n F = \e^{\Omega} F \, \e^{-\Omega}.
\]   

If we invert the $d \exp_{\Omega}$ operator given by (\ref{eq.4c}), we get from (\ref{ccv.4}) the formal identity
\[
\dot{\Omega} = \frac{x}{\e^x - 1} (A - \e^x F), 
\]
where $x \equiv \ad_{\Omega}$. Now, taking    
into account that
\[
\frac{x \, \e^x}{\e^x -1} = \frac{x}{\e^x - 1} + x,
\]   
it is clear that
\begin{equation} \label{eq.5a}
	\dot{\Omega} = \frac{x}{\e^x - 1} (A - F) - x F.
\end{equation}
A more convenient way of expressing this identity is obtained by recalling that
\[
d \exp_{\Omega}^{-1} (C) =  \frac{\ad_{\Omega} }{\exp(\ad_{\Omega}) - I} (C)= \sum_{k=0}^{\infty} \frac{B_k}{k!} \ad_{\Omega}^k (C),
\]
where $B_k$ are the Bernoulli numbers, so that (\ref{eq.5a}) reads
\begin{equation}  \label{eq.omega}
	\dot{\Omega} = \sum_{k=0}^{\infty} \frac{B_k}{k!} \ad_{\Omega}^k (A - F) \, - \,  \ad_{\Omega} F.
\end{equation}
With more generality, we can assume that both $A(t)$ and $F(t)$ are power series of some appropriate parameter $\varepsilon$,
\begin{equation}   \label{power.s1}
	A(t) = A_0(t) + \sum_{n \ge 1} \varepsilon^n A_n(t), \qquad\quad  F(t) = F_0(t) + \sum_{n \ge 1} \varepsilon^n F_n(t),
\end{equation}
and thus the generator $\Omega$ will be also a power series,
\begin{equation}  \label{power.s2}
	\Omega(t) = \sum_{n \ge 1}  \varepsilon^n \Omega_n(t), \qquad \Omega_n(0) = 0.
\end{equation}  
The successive $\Omega_n(t)$ can be determined by inserting (\ref{power.s1}) and (\ref{power.s2}) into equation (\ref{eq.omega}) and equating 
terms of the same power in $\varepsilon$. Thus, one obtains the following recursive procedure:

\begin{equation}  \label{recursion1}
	\begin{aligned}
		&  n = 0: \qquad F_0 = A_0  \\
		&  n = 1: \qquad \dot{\Omega}_1 + [\Omega_1, A_0] = W_1^{(0)}  \\
		&  n \ge 2: \qquad \dot{\Omega}_n + [\Omega_n, A_0] = W_n^{(0)} + G_n - V_n, 
	\end{aligned}
\end{equation}  

where
\begin{equation}   \label{recursion2}
	\begin{aligned}
		& W_n^{(0)} = A_n - F_n, \qquad n \ge 1  \\
		& W_n^{(k)} = \sum_{m=1}^{n-k} [ \Omega_m, W_{n-m}^{(k-1)} ], \qquad n \ge 1, \;\; 1 \le k \le n-1 \\
		& G_n = \sum_{k=1}^{n-1} \frac{B_k}{k!} W_n^{(k)}, \qquad V_n = \sum_{k=1}^{n-1} [ \Omega_k, F_{n-k}], \qquad n \ge 2 .  
	\end{aligned}
\end{equation}
Of course, although the change of variables $x \longmapsto X$ is completely general,  it is only worth to be considered 
if equation (\ref{ccv.2}) is simpler to solve than (\ref{ccv.1}). In the following we analyze
several ways of choosing $F$ fulfilling this basic requirement, and how in this way we are able to recover different exponential perturbative expansions.

%%%%%%%%%%%%%%%%%%%%%%%%%%%%%%%%%%%%%%%%%%
\section{Magnus expansion}

The simplest choice one can imagine of is taking $F = 0$ or, in other words, one is looking for a linear transformation rendering the original system (\ref{ccv.1}) into
\begin{equation}   \label{trivial}
	\frac{dX}{dt}=0,
\end{equation}
with trivial solution $X(t) = X(0) = x_0$. A sufficient condition for the reducibility of equation (\ref{ccv.1}) to (\ref{trivial}) is \cite{nemytskii89qto}
\[
\int_0^{+\infty} \|A(t)\|_F \, dt < +\infty,
\]
where $ \|A(t)\|_F = \sqrt{\sum_{i,j=1}^d |a_{ij}|^2}$.
If this is the case, from (\ref{ccv.3}),
\[
x(t) =   \exp(\Omega(t)) X(t) = \exp(\Omega(t)) x_0,
\]
where $\Omega(t)$ is determined from (\ref{eq.omega}) with $F = 0$, i.e.,
\begin{equation}   \label{ma.1}
	\dot{\Omega} = \sum_{k=0}^{\infty} \frac{B_k}{k!} \ad_{\Omega}^k A.
\end{equation}    
This, of course, corresponds to the well known Magnus expansion for the solution $x(t)$ of (\ref{ccv.1}) \cite{magnus54ote,blanes09tme}. 
The terms $\Omega_n(t)$ are then determined by 
the recursion (\ref{recursion2}) by taking $F_n=0$. If we take $A_0 = 0$ and $A_1(t) = A(t)$ in (\ref{power.s1}), then we get the familiar recursive
procedure \cite{blanes09tme}
\begin{equation}   \label{magnus1}
	\begin{aligned}
		& W_1^{(0)} = A,  \qquad W_n^{(0)} = 0, \qquad n \ge 2  \\
		& W_n^{(k)} = \sum_{m=1}^{n-k} [ \Omega_m, W_{n-m}^{(k-1)} ], \qquad n \ge 1, \;\; 1 \le k \le n-1 \\
		& G_n = \sum_{k=1}^{n-1} \frac{B_k}{k!} W_n^{(k)},  \qquad n \ge 2   \\
		&  \dot{\Omega}_1 = A_1, \qquad \dot{\Omega}_n = G_n, \qquad n \ge 2,
	\end{aligned}
\end{equation}
whence the successive terms $\Omega_n$ are obtained by integration. An explicit expression for $\Omega_n(t)$ involving only independent nested commutators 
can be obtained by working out this recurrence and using the class of bases proposed by  Dragt \& Forest \cite{dragt83con} for the
Lie algebra generated by
the operators $A_1(t_1), \ldots A_1(t_n)$. Specifically, one has \cite{arnal18agf}
\begin{equation}   \label{me.com3}
	\begin{aligned}
		& \Omega_n(t) = \frac{1}{n} \sum_{\sigma \in S_{n-1}} \, (-1)^{d_{\sigma}} \frac{1}{\binom{n-1}{d_{\sigma}}} \, 
		\int_0^t dt_1 \int_0^{t_1} dt_2 \cdots \int_0^{t_{n-1}} dt_n \,  \\
		&  \qquad\qquad\qquad  [A(t_{\sigma(1)}), [A(t_{\sigma(2)}) \cdots 
		[A(t_{\sigma(n-1)}), A(t_n)] \cdots ]],
	\end{aligned}    
\end{equation}
where $\sigma$ is a permutation of $\{1, 2, \ldots, n-1 \}$ and $d_{\sigma}$ corresponds to the number descents of $\sigma$. 
We recall that $\sigma$ has a descent in $i$ if $\sigma(i) > \sigma(i+1)$, $i=1,\ldots,n-2$. Notice that the argument in the last term is fixed to $t_n$, and
one considers all permutations in $A(t_1), A(t_2), \ldots, A(t_{n-1})$.
Moreover, 
the series (\ref{power.s2}) converges in this case in the interval $t \in [0, t_\mathrm{f}]$ such that
\begin{equation}   \label{conv.m}
	\int_0^{t_\mathrm{f}} \|A(s)\|_2 \, ds < \pi
\end{equation}
and the sum $\Omega(t)$ satisfies $\exp(\Omega(t)) = U(t)$ \cite{blanes09tme}. Here $\| \cdot \|_2$ denotes the spectral norm, 
characterized as
$\|A\|_2 = \max\{ \sqrt{\lambda}: \lambda \, \mbox{ is an eigenvalue of } \, A^{\dagger} A \}$.

The Magnus expansion has a long history as a tool to approximate solutions in a wide spectrum of fields in Physics and Chemistry, from atomic and molecular 
physics to Nuclear Magnetic Resonance to Quantum Electrodynamics (see \cite{blanes09tme} and references therein). Also in computational mathematics
it has been used to construct efficient algorithms for the numerical integration of differential equations within the widest field of geometric numerical integration
\cite{iserles00lgm,hairer06gni,blanes16aci}. Recently it has also been used in the treatment of dissipative driven two-level systems 
\cite{begzjav20mea}.

%%%%%%%%%%%%%%%%%%%%%%%%%%%%%%%%%%%%%%%%%%
\section{Floquet--Magnus expansion}
\label{FMsection}

The Magnus expansion can be in principle applied for any particular time dependence in $A(t)$, as long as the integrals in $\Omega_n(t)$ can be computed or
conveniently approximated. When $A(t)$
in (\ref{ccv.1}) is periodic in $t$ with period $T$, however, other changes of variables may be more suitable. According to Floquet's theorem, 
the original system is reducible to a system with a constant coefficient matrix $F$, whose eigenvalues (the so-called characteristic exponents) 
determine the asymptotic stability of the solution $x(t)$. In addition, the linear transformation is periodic with the same period $T$ \cite{hale80ode}.

In our general framework, then, it makes sense to determine
a change of variables 
$x =  \exp(\Omega(t)) X(t)$ in such a way that $F$ in (\ref{ccv.2}) is constant, so that $X(t) = \exp(t F) x_0$ and $\Omega(t)$ is periodic. Proceeding as before, if
we take $A_0 = 0$ and $A_1(t) = A(t)$ in (\ref{power.s1}), the procedure (\ref{recursion1})--(\ref{recursion2}) simplifies to
\begin{equation}   \label{floquet-magnus1}
	\begin{aligned}
		& W_1^{(0)} = A - F_1,  \qquad W_n^{(0)} = -F_n, \qquad n \ge 2  \\
		& W_n^{(k)} = \sum_{m=1}^{n-k} [ \Omega_m, W_{n-m}^{(k-1)} ], \qquad n \ge 1, \;\; 1 \le k \le n-1 \\
		& G_n = \sum_{k=1}^{n-1} \frac{B_k}{k!} W_n^{(k)},  \qquad V_n = \sum_{k=1}^{n-1} [ \Omega_k, F_{n-k}],   \qquad n \ge 2   
	\end{aligned}
\end{equation}
and
\begin{equation}  \label{floquet-magnus2}
	\begin{aligned}
		&  n = 1: \qquad \dot{\Omega}_1  = A - F_1 \\
		&  n \ge 2: \qquad \dot{\Omega}_n  =  - F_n + G_n - V_n \equiv \mathcal{F}_n - F_n
	\end{aligned}
\end{equation}  
Notice that, in general, $ \mathcal{F}_n = G_n-V_n$ depends on the previously computed $\Omega_k$, $F_k$, $k=1,\ldots,n-1$, 
so that equations (\ref{floquet-magnus2})
can be solved recursively as follows. First we determine $F_1$ and $F_n$ by taking the average of $A$ and $\mathcal{F}_n$, respectively, over one period $T$,
\[
F_1 =  \frac{1}{T} \int_0^T A(s) \, ds, \qquad      F_n = \langle \mathcal{F}_n \rangle \equiv \frac{1}{T} \int_0^T (G_n(s) - V_n(s)) \, ds,
\]
and then compute the integrals
\[
\Omega_1(t) = \int_0^t A(s) \, ds - t F_1, \qquad \Omega_n(t) = \int_0^t (G_n(s) - V_n(s)) \, ds - t F_n,
\]
respectively, thus ensuring that $\Omega_n$ is periodic with the same period $T$. 
This results in the well known Floquet--Magnus expansion  for the solution of  (\ref{ccv.1}),
\[
x(t) =     \exp(\Omega(t)) \, \exp(t F) x_0,
\]
originally introduced in \cite{casas01fte} and subsequently applied in different areas \cite{eckardt17aqg,haga19dot,mananga16otf,ota19feo}. 
In the context of periodic quantum mechanical problems, $H_{\mathrm{ef}} \equiv i
\hbar F$ is called the effective Hamiltonian of the problem. {This expansion presents the great advantage that, in addition to preserving unitarity as the
	Magnus expansion, also allows one to determine the stability of the system by analyzing the eigenvalues of $F$. }

As shown in  \cite{casas01fte}, the resulting series for $\Omega(t)$ is absolutely convergent
at least for $t \in [0, {t_\mathrm{f}}]$ such that
\[
\int_0^{t_\mathrm{f}} \|A(s)\|_2 \, ds < 0.20925...
\]    

The procedure can be easily generalized to quasi-periodic problems \cite{verdeny16qpd}. We recall that $A(t)$ is quasi-periodic in $t$ with frequencies $(\omega_1, \ldots, \omega_r)$ if
$A(t) = \tilde{A}(\theta_1, \ldots, \theta_r)$, where $\tilde{A}$ is $2\pi$-periodic with respect to $\theta_1, \ldots, \theta_r$ and $\theta_j = \omega_j t$ for
$j=1, \ldots, r$. In that case we can write
\[
A(t) = \sum_{k \in \mathbb{Z}^r} a_k \e^{i (k,\omega) t}
\]
where $(k, \omega) \equiv k_1 \omega_1 + \cdots + k_r \omega_r$ and $\sum |a_k|^2 < \infty$ \cite{corduneanu68apf}.

In this case $\mathcal{F}_n$ is also quasi-periodic (by induction)  and we take
\[
F_n =  \langle \mathcal{F}_n \rangle \equiv \lim_{T \rightarrow \infty} \frac{1}{T} \int_a^{a+T} \mathcal{F}_n(s) ds,
\]
the limiting mean value of $\mathcal{F}_n(t)$, independent of the particular value of $a$. In consequence,
\[
\dot{\Omega}_n = - F_n  + \mathcal{F}_n(t) =   \sum_{k \in \mathbb{Z}^r  \setminus \{0\} } f_{j,k} \, \e^{i (k,\omega) t}
\]   
and 
\[
\Omega_n(t) = \int_0^t \mathcal{F}_n(s) ds \, - t F_n
\]   
is also quasi-periodic with the same basic frequencies as $\mathcal{F}_n$.

%%%%%%%%%%%%%%%%%%%%%%%%%%%%%%%%%%%%%%%%%%
\section{Perturbed problems}

There are many problems that can be formulated as (\ref{ccv.1}) with $A(t)$ depending on a perturbation parameter $\varepsilon$,
\[
A(t) = A_0 + \sum_{n \ge 1} \varepsilon^n A_n(t),
\]
so that $A_0$ is a constant matrix and $A_n(t)$ are periodic or quasi-periodic functions of $t$. Here again the issue of reducibility
has received much attention along the years, and the problem consists in determining both the linear transformation $P(t)$ and the constant matrix $F$. In this section
we consider several possible ways to proceed, depending on the final matrix $F$ one is aiming for.

\subsection{Removing the perturbation}

One obvious approach is to take $F = A_0$. In other words, we try to construct a transformation rendering the original system into another one in which the perturbation is removed \cite{casas12utd}. The transformation fulfilling this requirement is therefore
\begin{equation}   \label{rm.1}
	x(t) = \exp(\Omega(t,\varepsilon)) \exp(t A_0) x_0, \qquad\qquad  \Omega(t,\varepsilon) = \sum_{n \ge 1} \varepsilon^n \, \Omega_n(t), 
\end{equation}
where we have written explicitly the dependence on $\varepsilon$, the parameter of the perturbation. 
The recurrence for determining the generators $\Omega_n$ is obtained from (\ref{recursion1})--(\ref{recursion2}) as
\begin{equation}   \label{rm.2}
	\begin{aligned}
		& W_n^{(0)} = A_n, \qquad n \ge 1  \\
		& W_n^{(k)} = \sum_{m=1}^{n-k} [ \Omega_m, W_{n-m}^{(k-1)} ], \qquad n \ge 1, \;\; 1 \le k \le n-1 \\
		& G_1 = 0, \qquad\qquad  G_n = \sum_{k=1}^{n-1} \frac{B_k}{k!} W_n^{(k)},  \qquad n \ge 2   \\
		&  \dot{\Omega}_n  + [\Omega_n, A_0]  = A_n + G_n, \qquad n \ge 1.
	\end{aligned}
\end{equation}
Alternatively, we can write
\begin{equation}  \label{eq.om1}
	\dot{\Omega}_n = \ad_{A_0} \Omega_n + A_n + G_n
\end{equation}
with solution verifying $\Omega_j(0) = 0$ given by
\[
\Omega_n(t) = \e^{t \, \ad_{A_0}} \int_0^t  \e^{-s \, \ad_{A_0}} \, \big(A_n(s) + G_n(s) \big) \, ds.
\]
As a matter of fact, this scheme can be related with the usual perturbation treatment of equation (\ref{ccv.1}). To keep the formalism  simple, let us
assume that $A(t)$ in (\ref{power.s1}) is $A(t) = A_0 + \varepsilon A_1(t)$ and work directly with the power series of the transformation,
\[
P(t,\varepsilon) = \exp( \Omega(t,\varepsilon) ) = \sum_{n \ge 0} \varepsilon^n P_n(t), \qquad \mbox{ with } \qquad P_0 = I.
\]
Then, by inserting $x(t) = P(t,\varepsilon) \, \e^{t A_0} x_0$ into equation (\ref{ccv.1}), one determines the differential equation satisfied by $P$ as
\[
\dot{P} + P A_0 = A P,
\]
whence the successive terms $P_n(t)$ verify
\[
\dot{P}_n + [P_n, A_0] = A_1 P_{n-1}, \qquad P_n(0)=0, \qquad n \ge 1.
\]
The solution, as with equation (\ref{eq.om1}), is given by
\[
P_n(t) = \e^{t A_0} \, g_n(t) \, \e^{-t A_0},
\]
with
\begin{equation}  \label{gs}
	\begin{aligned}
		g_n(t) & =\int_0^t ds \, \e^{-s A_0} A_1(s) P_{n-1}(s) \e^{s A_0} =   \int_0^t ds \, \e^{-s A_0} A_1(s) \e^{s A_0} g_{n-1}(s)      \\
		&  =       \int_0^t ds_1 \int_0^{s_1} ds_2 \ldots \int_0^{s_{n-1}} ds_n  \, A_I(s_1) A_I(s_2) \ldots A_I(s_n).
	\end{aligned}
\end{equation}
Here $g_0 \equiv I$ and we have denoted $A_I(s) \equiv \e^{-s A_0} A_1(s) \, \e^{s A_0}$. Therefore, the solution of (\ref{ccv.1}) reads
\[
\begin{aligned}
x(t) & = P(t,\varepsilon) \, \e^{t A_0} x_0 = \left( I + \sum_{n \ge 1} \varepsilon^n P_n(t) \right)  \e^{t A_0} x_0   =    
\left( \e^{t A_0} + \sum_{n \ge 1} \varepsilon^n \e^{t A_0} g_n(t) \right) x_0 \\
& = \e^{t A_0} \left( I + \sum_{n \ge 1} \varepsilon^n g_n(t) \right) x_0
\end{aligned}     
\]
and so, taking into account the explicit expression (\ref{gs}) for $g_n(t)$, one recovers the solution provided by the standard perturbation theory in the
picture defined by $A_0$.

\subsection{Lie--Deprit perturbative algorithm}

The previous treatment has one important drawback (in addition to the lack of preservation of unitarity when one deals with the
Schr\"odinger equation and the expansion is truncted), namely, when $A(t)$ is a periodic or quasi-periodic function of time, the secular terms are
not removed and, as a result, the time interval of validity of the resulting approximations 
is very small indeed. For this reason it is worth considering in the general framework
(\ref{recursion1})--(\ref{recursion2}) the quasi-periodic case and look for a quasi-periodic transformation leading to a \textit{constant} coefficient matrix
\begin{equation}   \label{Fcte}
	F(\varepsilon) =  F_0 + \sum_{n \ge 1} \varepsilon^n F_n, \qquad \mbox{ with } \qquad F_0 = A_0.
\end{equation}
This problem has been addressed a number of times in the literature (see e.g. \cite{burd07moa,jorba92otr} and references therein). In particular, in 
\cite{arnal16apa} a perturbative algorithm is presented for constructing this  transformation as the exponential of a quasi-periodic matrix. Here we show
how the results of \cite{arnal16apa,casas12utd} can be reproduced by the generic scheme proposed here in a more direct way, in the sense that the terms
$\Omega_j$ are determined at once from (\ref{recursion1})--(\ref{recursion2}).

It is clear that the problem reduces to solve
\[
\dot{\Omega}_n + [\Omega_n, A_0] = \mathcal{F}_n - F_n, \qquad \Omega_n(0) = 0
\]
where now
\[
\mathcal{F}_1 \equiv  A_1,  \qquad\qquad \mathcal{F}_n \equiv   A_n + G_n - V_n, \quad n \ge 2,
\]   
and the goal is to determine the constant term $F_n$ and construct $\Omega_n(t)$ as a quasi-periodic function with
the same basic frequencies as $A_n$ for $n \ge 1$. 

As shown in the Appendix, the solution verifying all these requirements is given by
\[
F_n = \langle \mathcal{F}_n(t) \rangle - [A_0, M_n(0)], \qquad\quad  \Omega_n(t) = - M_n(0) + \e^{t \, \ad_{A_0}} \, M_n(t),
\]
where $M_n(t)$ is the antiderivative
\[
M_n(t) = \int   \e^{-t \, \ad_{A_0}} \big( \mathcal{F}_n(t) - \langle \mathcal{F}_n(t) \rangle \big) dt.
\]   
Finally, the solution of (\ref{ccv.1}) is written as
\[
x(t) = \exp(\Omega(t,\varepsilon)) \exp(t F(\varepsilon)) x_0,
\]
where $\Omega(t,\varepsilon)$ and $F(\varepsilon)$ are appropriate truncations of the corresponding series.

\subsection{Quantum averaging}
In the context of time-dependent quantum mechanical systems, the averaging method has been used to construct quantum
analogues of classical perturbative treatments \cite{scherer96qa,scherer97npa}. Essentially, the basic approach in that setting is to transform the original
Hamiltonian of the problem by a unitary transformation so that the problem of finding the time evolution of the
transformed Hamiltonian is easier than the original one. The idea has been also applied in the perturbative treatment
of pulse-driven quantum problems \cite{daems03otd,daems04pdq,sugny04tdu}.

This approach to quantum averaging can be fit into our general framework of Section \ref{sec.2} by taking $F_0 = A_0$ and the terms $F_n(t)$ in (\ref{power.s1})
verifying in addition 
\begin{equation}  \label{condi.f}
	\dot{F}_n(t) + [F_n(t), A_0] = 0, \qquad n \ge 1.
\end{equation}
Clearly, the solution of (\ref{condi.f}) verifies $F_n(t) = \e^{t \ad_{A_0}} F_n(0) = \e^{t A_0} F_n(0) \e^{-t A_0}$, and thus
\[
\frac{d }{dt} \left( \e^{-t A_0} X(t) \right) = \left( F(0) - A_0 \right) \, \e^{-t A_0} X(t),
\]
which leads to $\exp(-t A_0) X(t) = \exp(t (F(0) - A_0) X(0)$ and finally
\[
X(t) = \e^{t A_0} \e^{t (F(0) - A_0)} X(0).
\]
In other words, condition (\ref{condi.f}) guarantees that equation (\ref{ccv.2}) can be solved in a closed way and the dynamics of (\ref{ccv.1}) is obtained once
the generator $\Omega(t)$ is determined, even when $F_n$ depends explicitly on time.  As shown in \cite{scherer97npa}, the corresponding solutions are
\[
F_n(t)  = \lim_{T \rightarrow \infty} \frac{1}{T} \int_0^T B_n(\sigma,t) d \sigma, \qquad 
\Omega_n(t)  = \lim_{T \rightarrow \infty} \frac{1}{T} \int_0^T d \tau \int_0^{\tau} d\sigma \big( B_n(\sigma,t) - F_n(t) \big),
\]
where 

\[
B_n(\sigma,t) = \e^{(t- \sigma) A_0} \mathcal{F}_n(t-\sigma)   \e^{-(t- \sigma) A_0}
\]
provided 

\[
\lim_{T \rightarrow \infty}   \frac{1}{T} \big( \mathcal{F}_n(t) - \e^{-T A_0} \mathcal{F}_n(t-T) \e^{T A_0} \big) = 0.
\]      
Here $\mathcal{F}_1 =  A_1$ and $\mathcal{F}_n =   A_n + G_n - V_n$ for  $n \ge 2$.

%%%%%%%%%%%%%%%%%%%%%%%%%%%%%%%%%%%%%%%%%%

\section{Illustrative examples}

In previous sections we have reviewed several perturbative expansions aimed to get analytic approximations for the solution of (\ref{ecu.1}), and how they 
can be derived from a same basic principle, namely a transformation of coordinates. This allows one, in particular, to design a unique 
computational procedure to deal with a particular problem defined by $A(t)$ and take the most appropriate variant in each case. To illustrate the technique
we next consider two examples describing the dynamics of simple time-dependent quantum systems, although the same treatment can be applied of course
to other more involved problems.

\subsection{The three-lambda system}

The so-called driven three-lambda system describes an atomic three energy-level system with two ground states $|1\rangle$, $|2 \rangle$ 
with the same energy $E_1$ that are coupled
with an excited state $| 3 \rangle$ with energy $E_3$ via a time-dependent laser field.  In the interaction picture the Hamiltonian is given by
$H(t) = f(t) \, |3\rangle \big( \langle 1| + \langle 2 | \big)$ (+ Hermitian conjugate), or equivalently,
\begin{equation}  \label{3lambda.11}
	H(t) =  \left(  \begin{array}{ccc}
		0  &  0  &  \bar{f}(t)  \\
		0  &  0  &   \bar{f}(t)  \\
		f(t)  &  f(t)  &  0
	\end{array} \right),	
\end{equation}
where $f(t)$ is usually a periodic function. The corresponding Schr\"odinger equation, $i \hbar \dot{U}(t) = H(t) U(t)$, is then a particular case of equation (\ref{ccv.1})
with $A(t) = -i H(t)/\hbar$, and one is typically interested in obtaining the induced probability of transition between states $|1\rangle$ and $|2 \rangle$, namely
$P_{12}(t) = |\langle 1 |U(t)| 2 \rangle|^2$.

\begin{figure}[H]
	\includegraphics[scale=0.3]{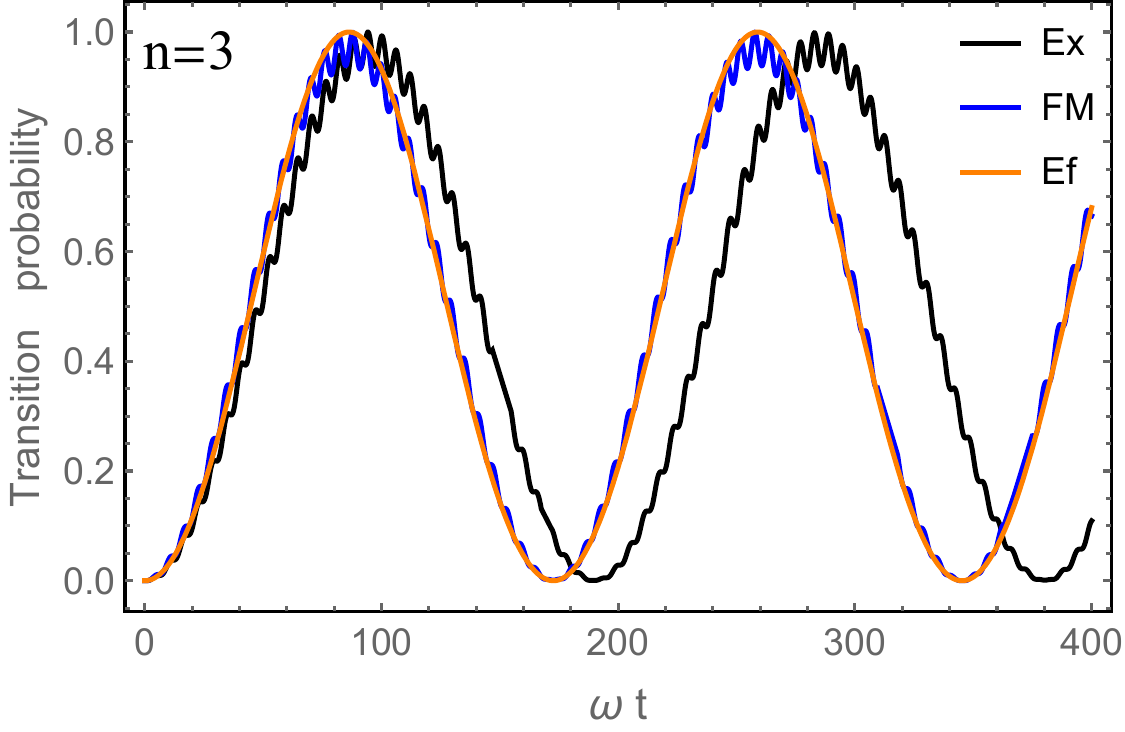}\includegraphics[scale=0.3]{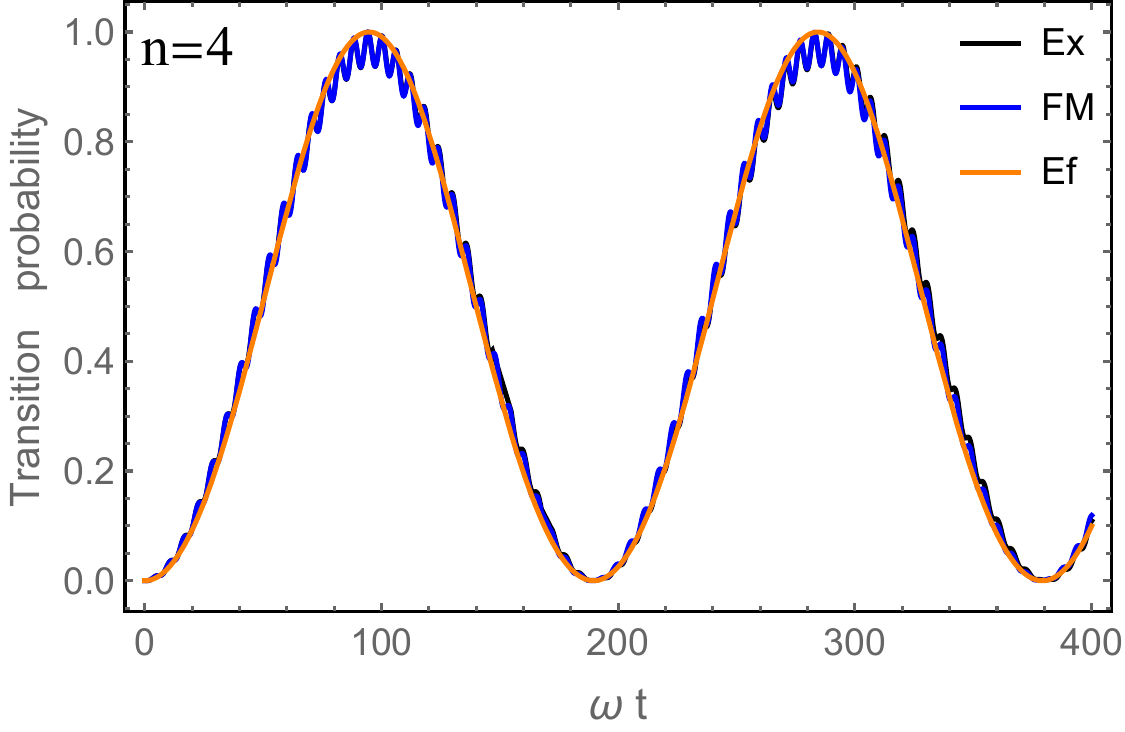}
	\caption{{\small Transition probability for the three-lambda system with $f(t) = \e^{i \omega t}$  for
			$\omega=10 \, (1+\frac{\sqrt{2}}{2})^{-\frac{1}{2}}$} obtained with Floquet--Magnus (FM) 
		with $n=3$ (left) and $n=4$
		(right) terms in the expansion. Blue line corresponds to FM,
		orange (Ef) line stands for the effective Hamiltonian $H_{\mathrm{ef}} \equiv i F$ and black (Ex) line
		is the exact solution.}	 	\label{figu11}	
\end{figure}

We follow \cite{verdeny16qpd} and take 
\begin{equation}   \label{3lambda.22}
	f(t)=\beta \, \e^{i \omega t}
\end{equation}    
with $\beta=1$ and $\omega=\frac{10}{(1+\frac{\sqrt{2}}{2})^\frac{1}{2}} \approx 7.65$ as an example of periodic function. Starting from the initial condition
$U(0) = I$ we compute the approximate solution until the final time $\omega \, t = 400$, i.e for 487 periods, with the Floquet--Magnus (FM) expansion 
(\ref{floquet-magnus1})--(\ref{floquet-magnus2}) up to $n=3$ and $n=4$ and determine the corresponding transition probability. 
In this way we get Figure \ref{figu11}, where the result achieved by the effective Hamiltonian 
$X(t) = \e^{t F}$, with $F = -i H_{\mathrm{ef}}/\hbar$, is also depicted. 
We see that for this value of $\omega$, considering more terms in the series for both $\Omega(t)$ and $F$ leads to better results and that 
working only with the effective Hamiltonian gives indeed a good approximation. The reference solution is computed with the \texttt{DSolve} function of 
\textit{Mathematica}. For completeness, the value of the effective Hamiltonian  up to $n=4$ reads 
\[
H_{\mathrm{ef}} = \left(
\begin{array}{ccc}
\frac{6-\omega ^2}{\omega ^4} & \frac{6-\omega ^2}{\omega ^4} & \frac{4}{\omega ^3} \\
\frac{6-\omega ^2}{\omega ^4} & \frac{6-\omega ^2}{\omega ^4} & \frac{4}{\omega ^3} \\
\frac{4}{\omega ^3} & \frac{4}{\omega ^3} & \frac{2 \left(\omega ^2-6\right)}{\omega ^4} \\
\end{array}
\right),
\]
whereas the transition probability obtained with $H_{\mathrm{ef}}$ up to $n=3$ is given by ($\hbar = 1$)
\[
P_{12}^{(3)}  = \frac{  \sin ^2(t \, \widetilde{\omega}) }{ (\omega^2 +8)}   \left(-4 \cos(2 t  \widetilde{\omega}) +\omega ^2+4\right)
\]   
with $\widetilde{\omega} \equiv \left(\sqrt{\omega^2+8} \right)/\omega ^3$.

\begin{figure}[H]
	\includegraphics[scale=0.35]{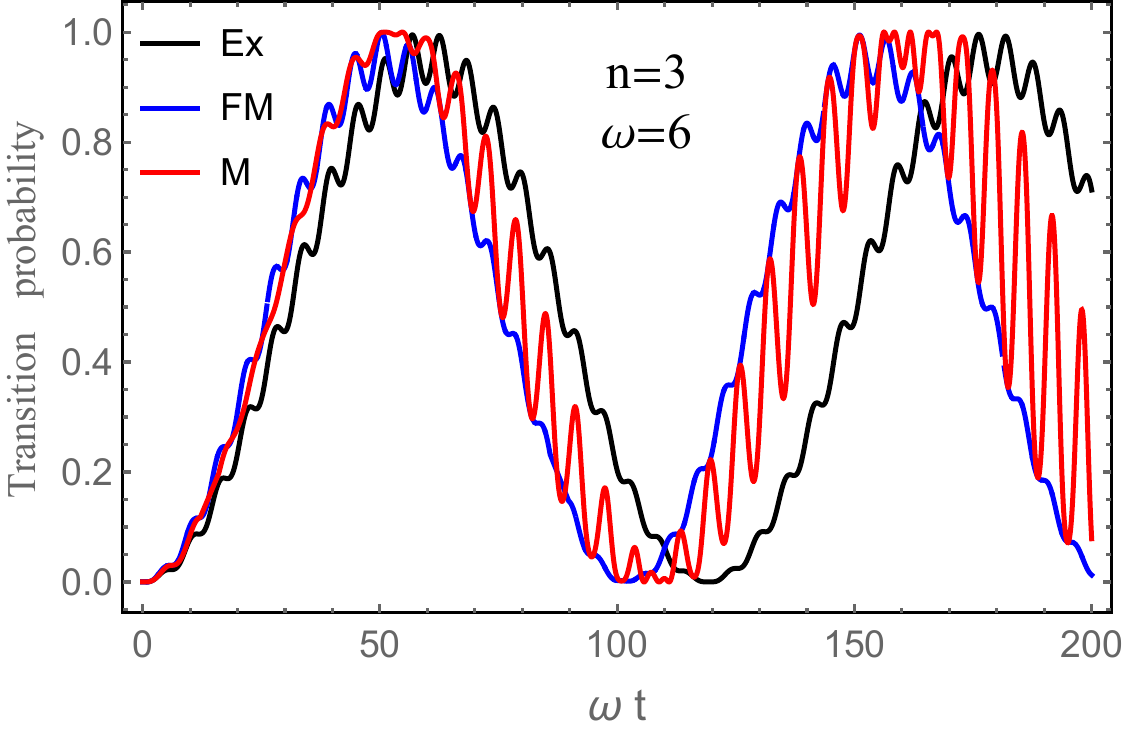}	\includegraphics[scale=0.35]{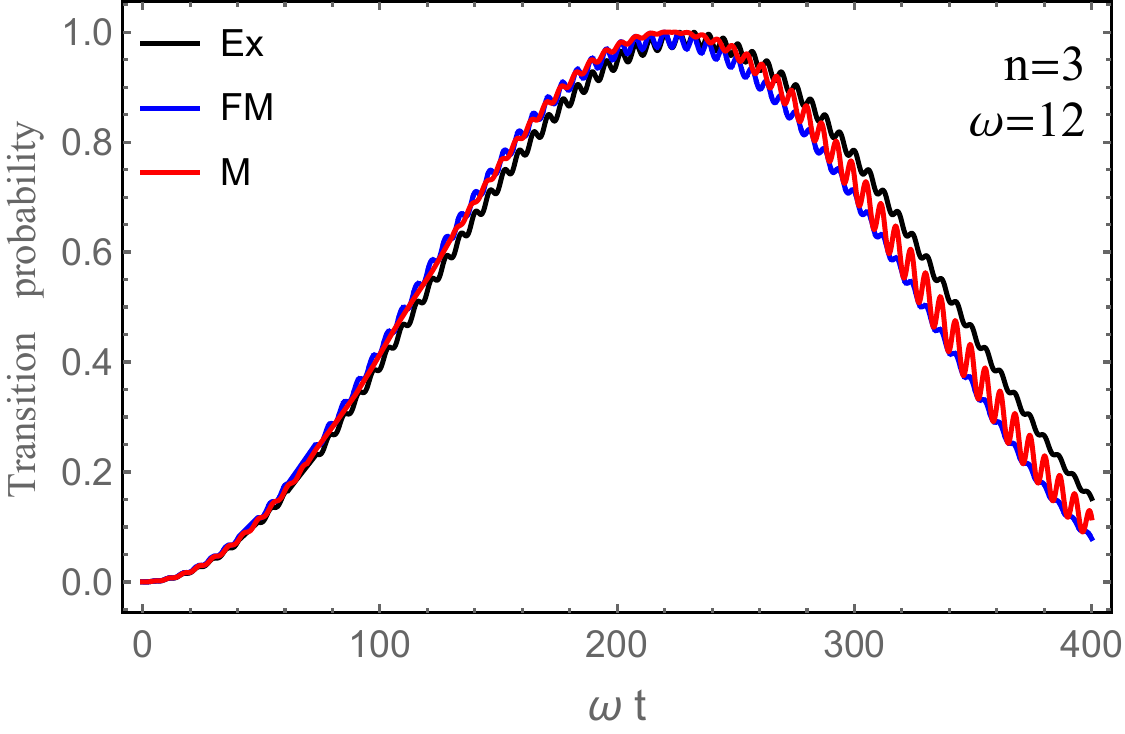}\\
	\includegraphics[scale=0.35]{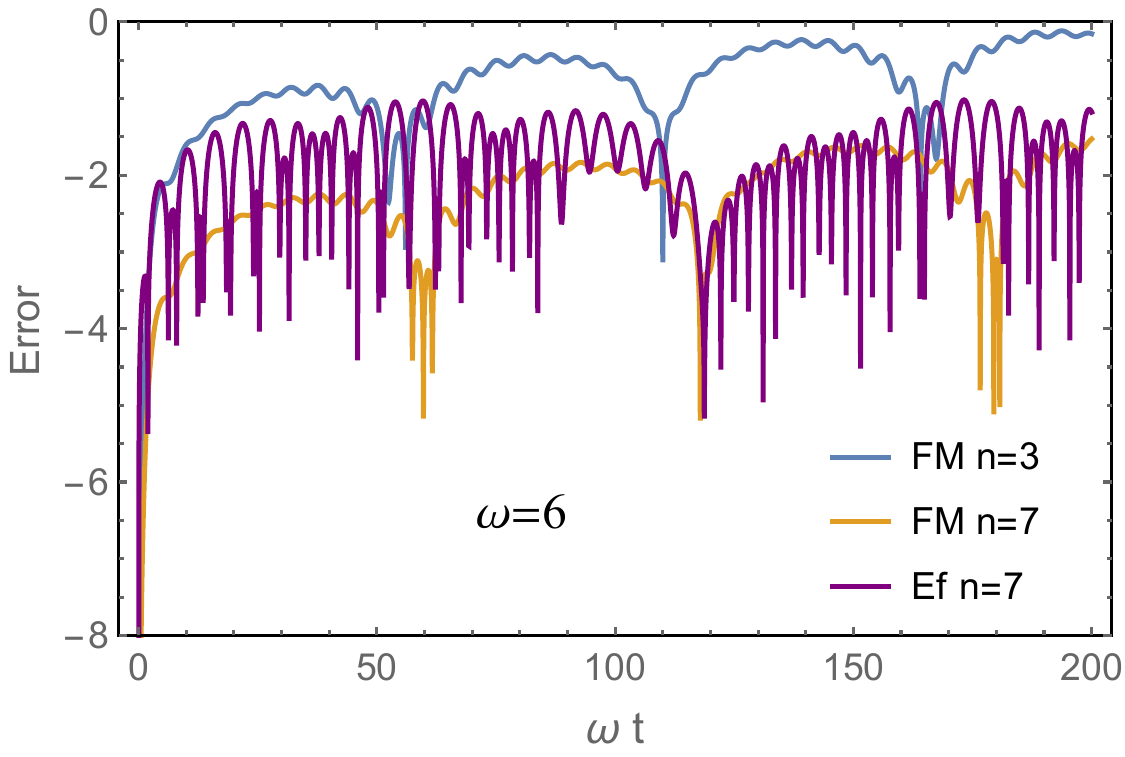} \includegraphics[scale=0.35]{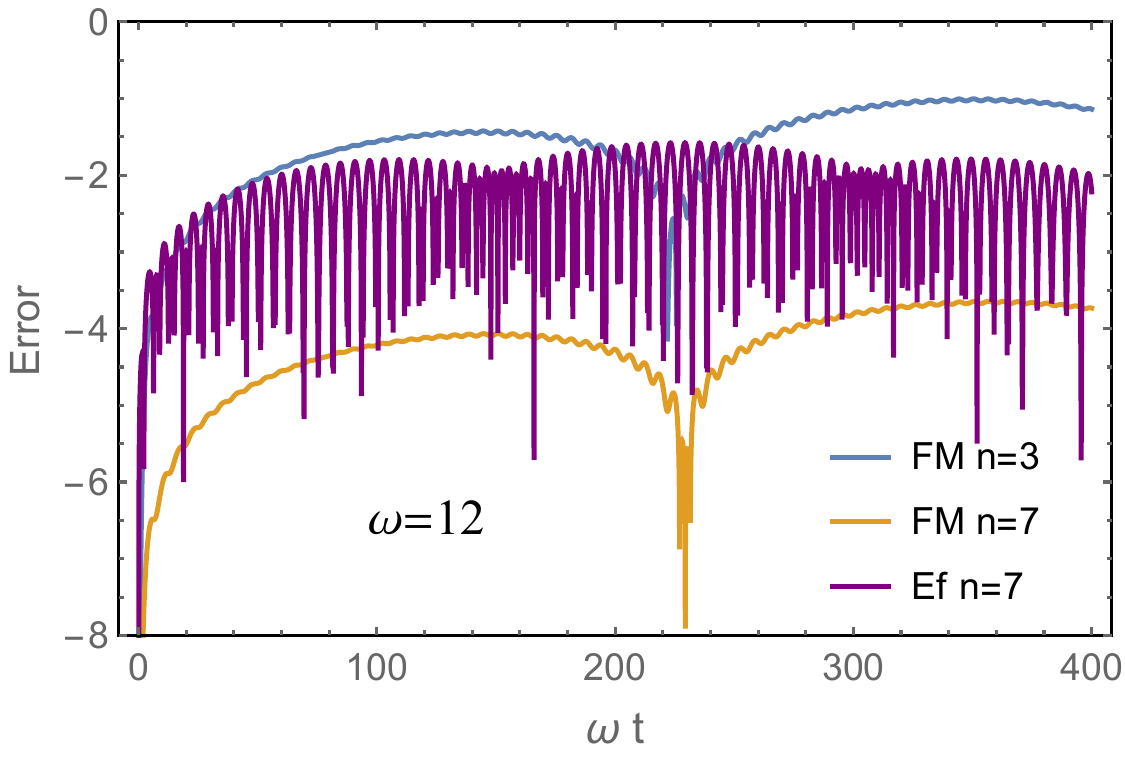}	\\
	\caption{{\small Top: transition probability for the three-lambda system with $f(t) = \e^{i \omega t}$  obtained with Floquet--Magnus (blue line) and Magnus 
			expansion (red line), both with $n=3$ terms. Left diagram is for $\omega = 6$, right diagram is for $\omega=12$. Bottom: error (in logarithmic scale)
			as a function of time
			for Floquet--Magnus with $n=3$ and $n=7$ terms, and the effective solution with $n=7$ terms (purple).}}
	\label{figu22}
\end{figure}

Next we repeat the simulation  with $\omega = 6$ and $\omega = 12$, by taking $n=3$ terms in both the Floquet--Magnus and the Magnus expansion
(\ref{magnus1}). The results are shown in Figure \ref{figu22} (top). Notice that the approximations get worse for smaller values of $\omega$. 
To get a more
quantitative view, also in Figure \ref{figu22} (bottom) we show the error in the transition probability (in logarithmic scale) for the same values of $\omega$, but
now include $n=3$ and $n=7$ terms in the FM expansion, as well as the result obtained with just the effective Hamiltonian in this last case.   As is clear
from the figure, taking into account more terms in the FM expansion improves the approximations, and this improvement is more remarkable for
larger values of $\omega$.

Here the approximate solution is obtained as 
\[
U(t) = \exp(\tilde{\Omega}(t)) \, \exp(-i \, t H_{\mathrm{ef}}/\hbar)
\]
and the two parts of $U$
encode different behaviors of the system. Thus, $\exp(\tilde{\Omega}(t))$ describes fast fluctuations around the envelope evolution provided by
$\exp(-i t H_{\mathrm{ef}}/\hbar)$, but these fluctuations are not always negligible, especially if the driving frequency is not too large.  In any case, it is worth noticing
that taking only into consideration the effective Hamiltonian renders less accurate results than working with the full approximation.

\begin{figure}
\includegraphics[scale=0.35]{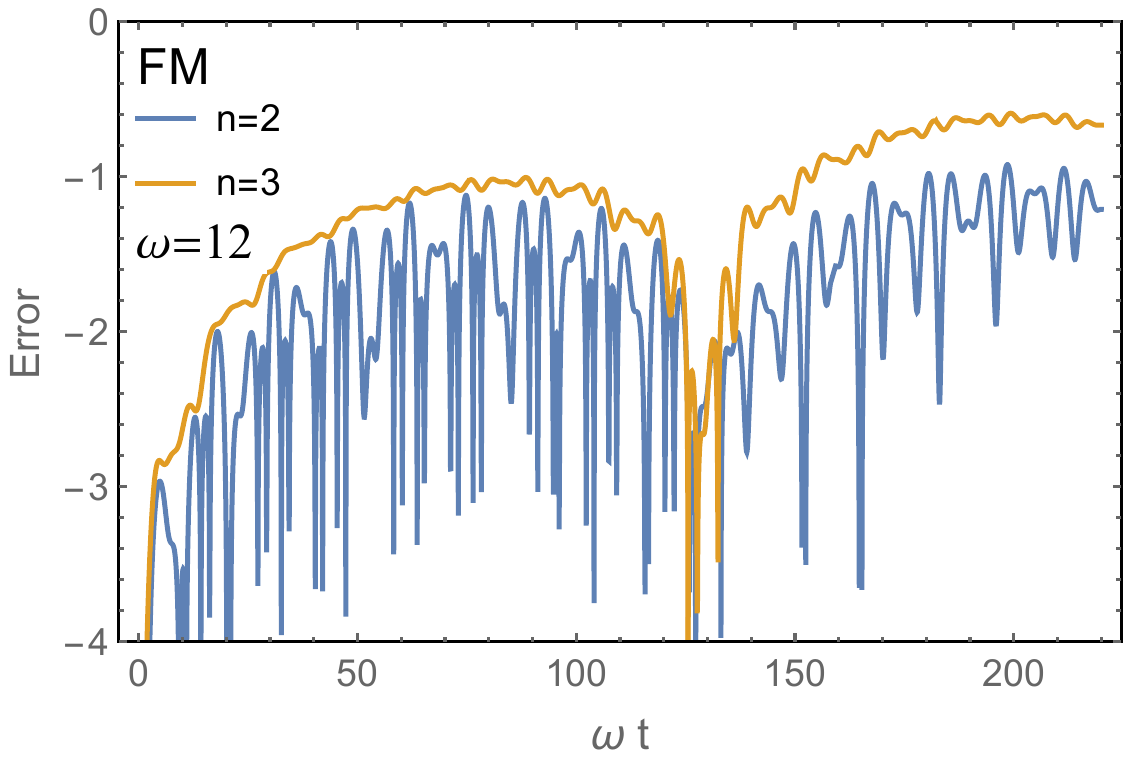} \includegraphics[scale=0.35]{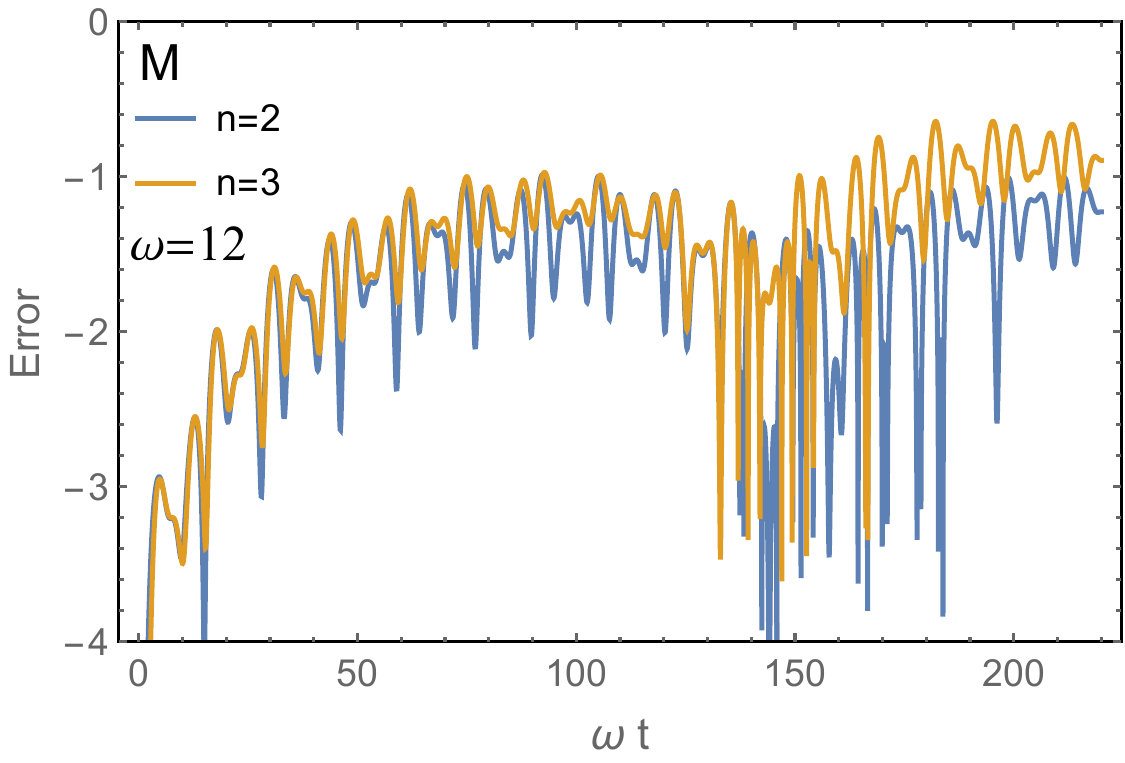}
\caption{{\small Error (in logarithmic scale) the transition probability (in logarithmic scale) for the quasi-periodic case (\ref{3lambda.33}) as a function of time
		for FM (left) and Magnus (right) with $n=2$ and $n=3$ terms in the expansion. In both cases $\omega = 12$ in (\ref{3lambda.33}).}} \label{figu33}
\end{figure}

\begin{figure}[H]
\begin{center}
	\includegraphics[scale=0.6]{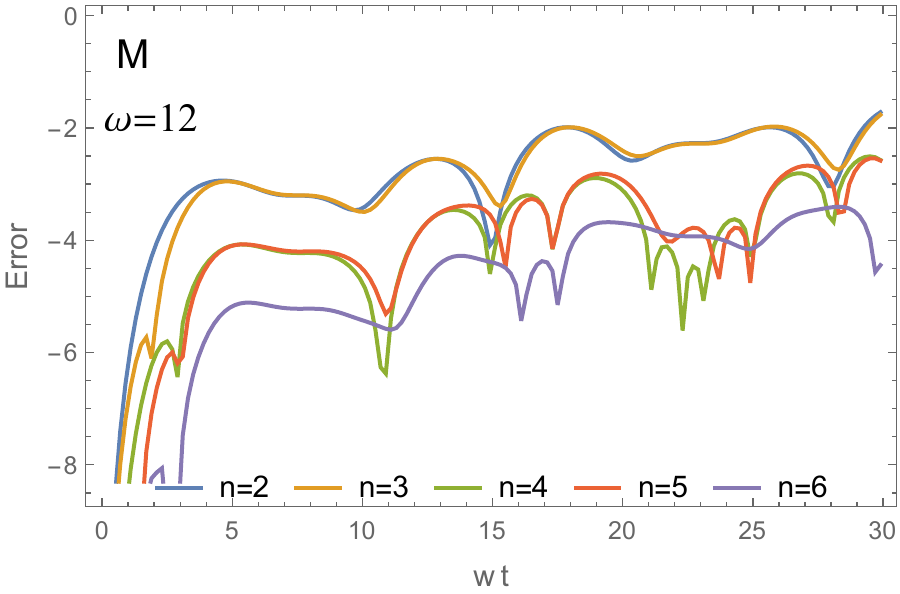} 
	\caption{{\small Detail of  Figure \ref{figu33} right in the interval $\omega t \in [0, 30]$. We have included up to $n=6$ terms in the Magnus
			series.}} \label{figu33c}
\end{center}	
\end{figure}

Next we take $f(t)$ in (\ref{3lambda.11}) to be a quasi-periodic function, namely
\begin{equation}  \label{3lambda.33}
	f(t) = \beta (   \e^{i \omega t} + \e^{i \sqrt{2} \omega t} ),
\end{equation}
with $\beta = 1$ and $\omega = 12$, and compute approximations using the Floquet--Magnus and Magnus expansions with $n=2$ and $n=3$ terms. 
The corresponding errors in the transition probability are shown in Figure \ref{figu33}. We see that in this case including more terms
in the expansions does not improve significantly the accuracy of the approximations. This is an indication of the lack of convergence of the expansions
for these values of the parameters. In fact, it is straightforward to verify that $\|A(t)\|_2 = \sqrt{8} \, \beta \left| \cos \left( \frac{\sqrt{2}-1}{2} \omega t
\right) \right|$, so that for the values of the parameters considered, the convergence of the Magnus expansion is ensured 
for $t < 1.6117$ according to the estimate (\ref{conv.m}). To illustrate the behavior in this range of times, in Figure \ref{figu33c} we collect the
results obtained by taking up to $n=6$ terms in the expansion. It is worth noticing the systematic decrease in the error when more terms are
included in the expansion. A similar pattern is observed for the Floquet--Magnus expansion, but with much smaller times (the convergence is
only assured for $t < 0.074$).

\subsection{Bloch--Siegert dynamics}

%In the Schr\"odinger picture, to illustrate perturbed problems with this approach

This is an example of a periodically driven two-level quantum system described by the time-dependent Hamiltonian
\begin{equation}   \label{bs.1}
H(t) = \left(  \begin{array}{cc}
\omega_0/2  &  2 b \cos(\omega t) \\
2 b \cos(\omega t) &  -\omega_0/2
\end{array} \right).
\end{equation}
The coupling parameter $b$ is the amplitude of a driving radio-frequency field and $\omega$ is its frequency. We can write
\[
H(t) = H_0 + \varepsilon H_1(t)
\]
with $\varepsilon = 2 b$ and 
\[
H_0 = \frac{\omega_0}{2} \sigma_3, \qquad \quad H_1(t) =  \cos(\omega t) \, \sigma_1,  
\]
where $\sigma_i$ are Pauli matrices.
With this simple example we can illustrate the different treatments previously considered, both perturbative and non-perturbative, and also analyze
how the results vary with the perturbation parameter $\varepsilon$.

\begin{figure}[H]
	\includegraphics[scale=0.4]{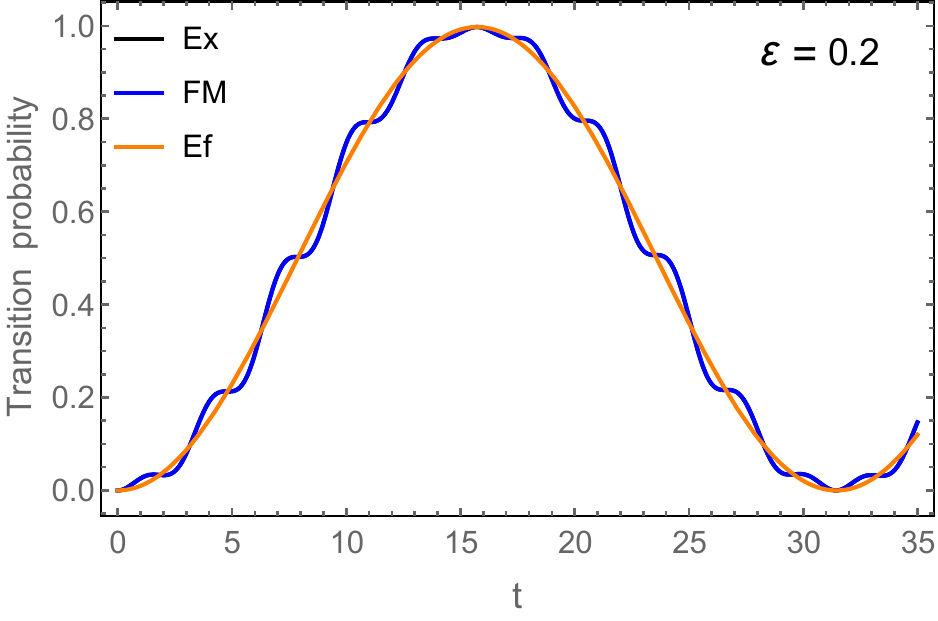}\includegraphics[scale=0.4]{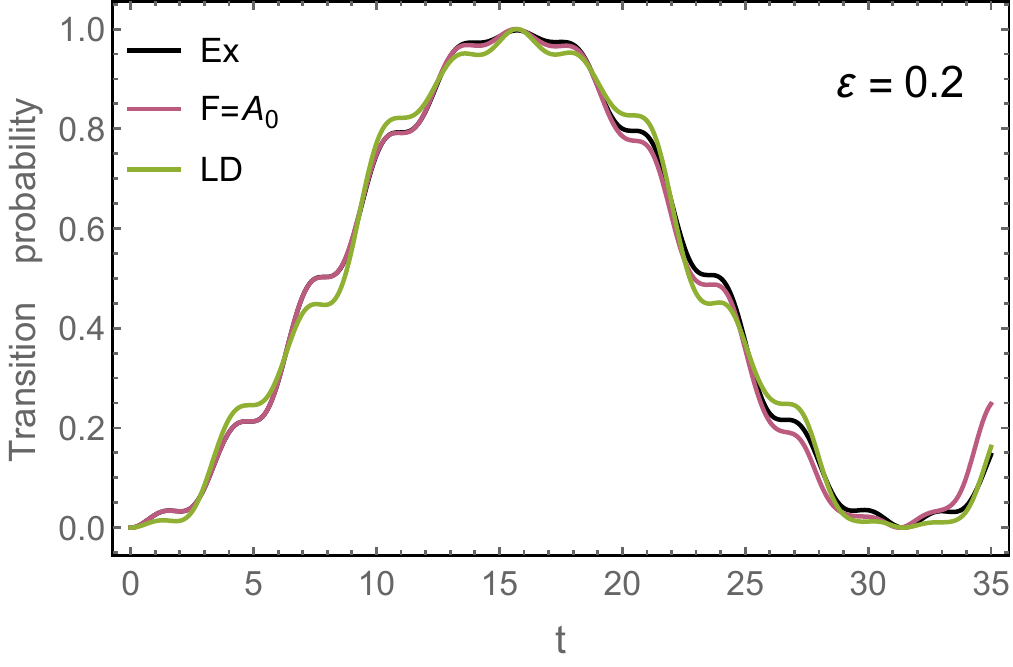}\\
	\includegraphics[scale=0.4]{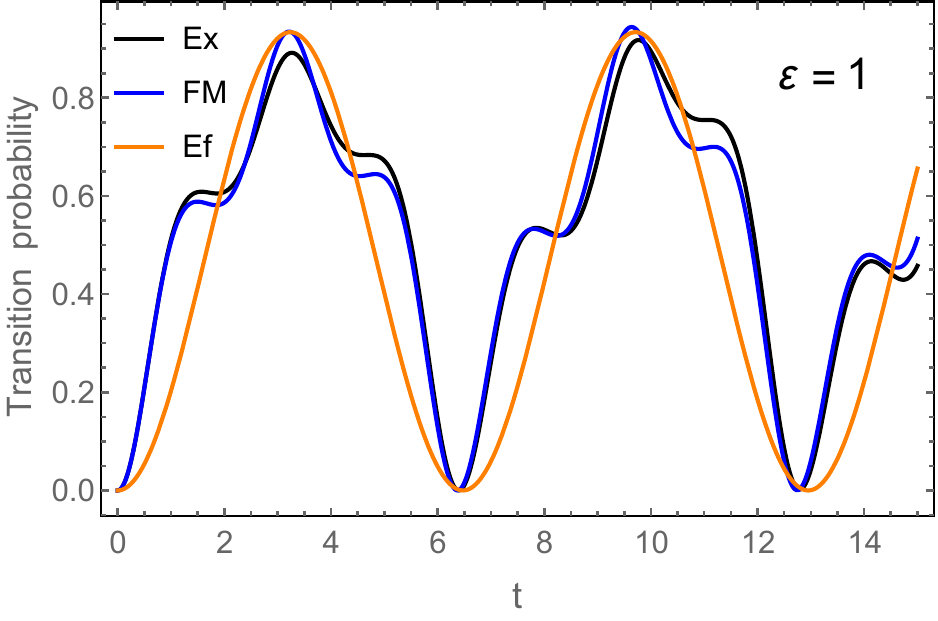}\includegraphics[scale=0.4]{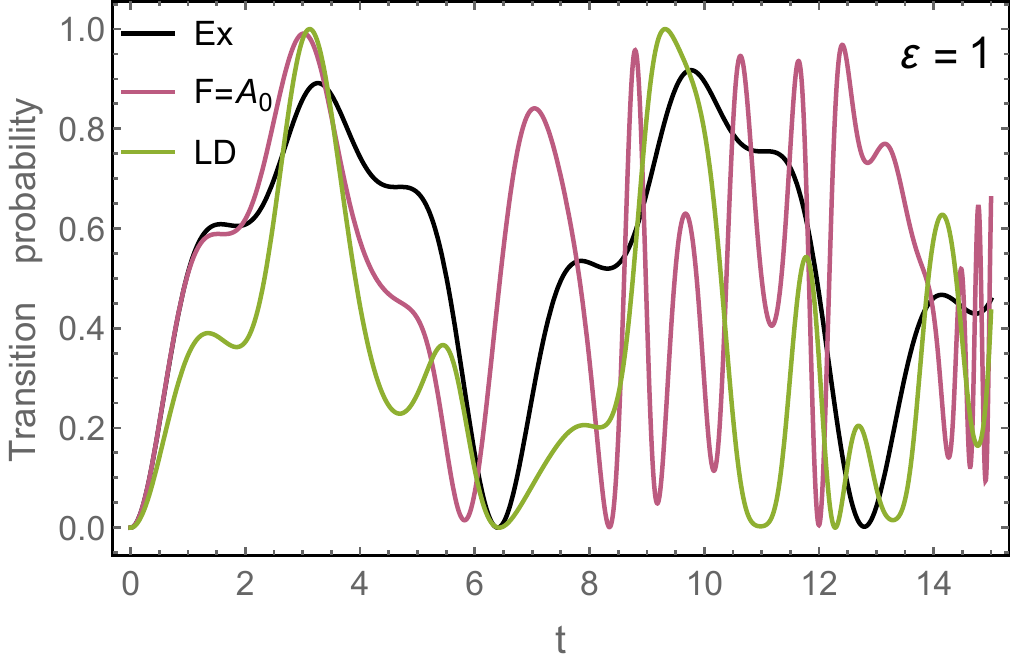}
	\caption{{\small Transition probability obtained with FM (left, blue line) and LD (right, green line) with $n=3$ terms in comparison with the exact result for the Bloch--Siegert Hamiltonian with $\varepsilon = 0.2$ (top) and $\varepsilon = 1$ (bottom). Black line corresponds to the exact solution, and orange to the
			result achieved with the effective Hamiltonian. On the right panels the results achieved by the algorithm removing the perturbation are also depicted
			($F = A_0$).}}  \label{figu44}
\end{figure}

As in \cite{giscard20gsf}, we consider the resonant case $\omega=\omega_0=1$ and two different values of the perturbation,
$\varepsilon = 0.2$ and $\varepsilon = 1$. One could use the Magnus expansion (\ref{magnus1}) directly to the matrix $A(t) = -i H(t)/\hbar$. In that case,
however, the convergence of the series, as given by (\ref{conv.m}), is only ensured for times $t \in [0, t_{\mathrm{f}}]$ such that 
\[
\int_0^{t_{\mathrm{f}}} \sqrt{ \frac{1}{4} + \varepsilon^2 \cos^2(t)} \, dt < \pi.
\]
Thus, in particular, if $\varepsilon = 0.2$, then $t_{\mathrm{f}} \approx 6.056$, whereas for $\varepsilon = 1$ it reduces to $t_{\mathrm{f}} \approx 3.608$.
We cannot expect, then, accurate results for larger time intervals, and this is indeed what is observed in practice. Something analogous happens also for the
Floquet--Magnus expansion, with the difference that the interval of convergence is even smaller. It makes sense, in consequence, to apply both procedures
in the interaction picture, i.e., to the matrix $A_I(t) = \e^{-t A_0} A_1(t) \e^{t A_0}$ and this is what we have done in our computations
by applying the procedure of section \ref{FMsection}.
We also compare with 
the Lie--Deprit (LD) perturbation algorithm and the formalism removing the perturbation, i.e., such that $F=A_0$.
As in the previous example, we compute $|U_{12}(t)|^2$, i.e., the probability of transition between the two levels, obtained with
each procedure with $n=3$ terms. The corresponding results are shown in Figure \ref{figu44} as a function of time, where we also include the exact result for comparison. 
Left diagrams are obtained with FM (in the interaction picture), 
whereas we depict in the right diagrams the results obtained by LD and also when $F$ is taken as $A_0$. Notice that  
top and bottom panels correspond to $\varepsilon = 0.2$ and
$\varepsilon = 1$, respectively. We have also depicted the result obtained with just the effective Hamiltonian in the
case of the Floquet--Magnus expansion. 

It can be observed that both the LD perturbative algorithm and the one obtained by removing the perturbation provide 
accurate results only when $\varepsilon$ takes small values, in contrast with  Floquet--Magnus in the interaction picture, with an enlarged validity domain in
$\varepsilon$. Notice
that also in this case the effective Hamiltonian is not enough to get a precise description of the dynamics when $\varepsilon$ increases. We should remark that
all schemes render unitary approximations by construction.

\begin{figure}[H]
	\includegraphics[scale=0.43]{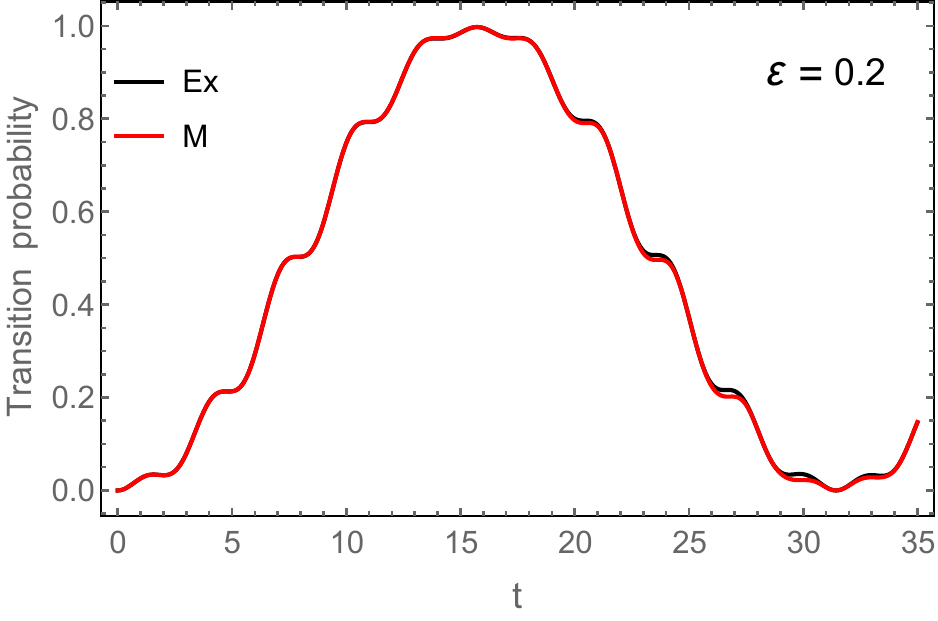}\includegraphics[scale=0.43]{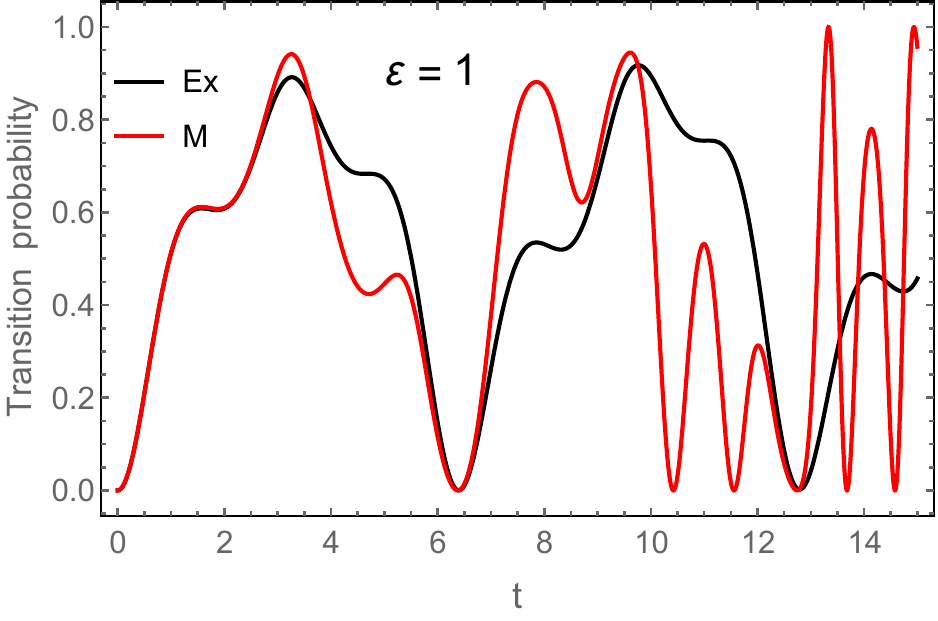}
	\caption{{\small Transition probability obtained with the standard Magnus expansion with $n=3$ terms (red line) in comparison with 
			the exact result (black line) for the Bloch--Siegert Hamiltonian with $\varepsilon = 0.2$ (left) and $\varepsilon = 1$ (right).}} \label{figu55}
\end{figure}

For comparison, in Figure \ref{figu55} we collect the results achieved by the standard Magnus expansion in the interaction picture
for this same problem and values of the perturbation
parameter: $\varepsilon = 0.2$ (left) and $\varepsilon = 1$ (right). Here again increasing the value of $\varepsilon$ leads to a loss of accuracy for long times.
If the Magnus expansion is applied to the matrix $A(t)$, then only accurate results are obtained for times in the convergence domain.

\begin{figure}[H]
	\includegraphics[scale=0.38]{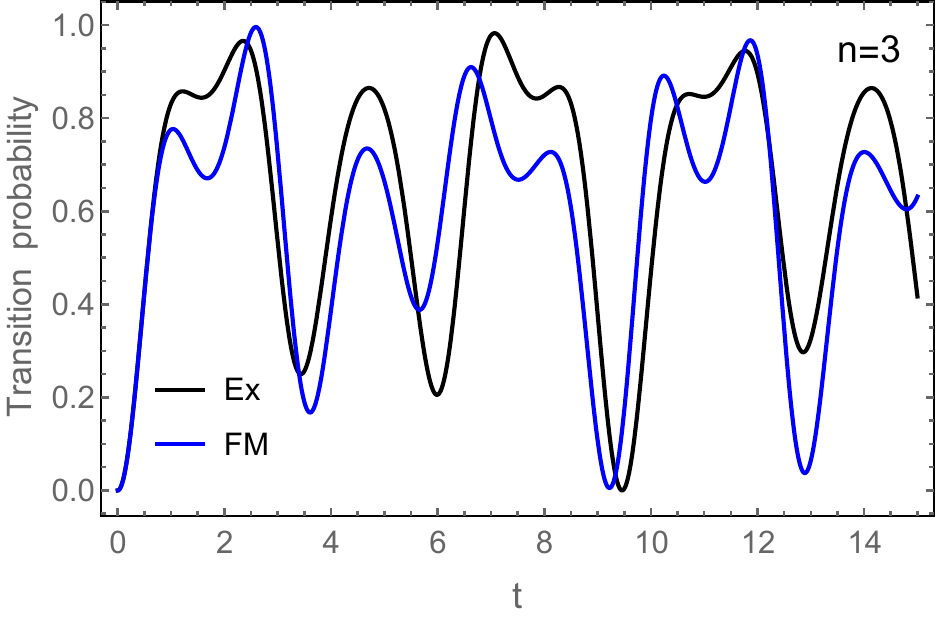}\includegraphics[scale=0.38]{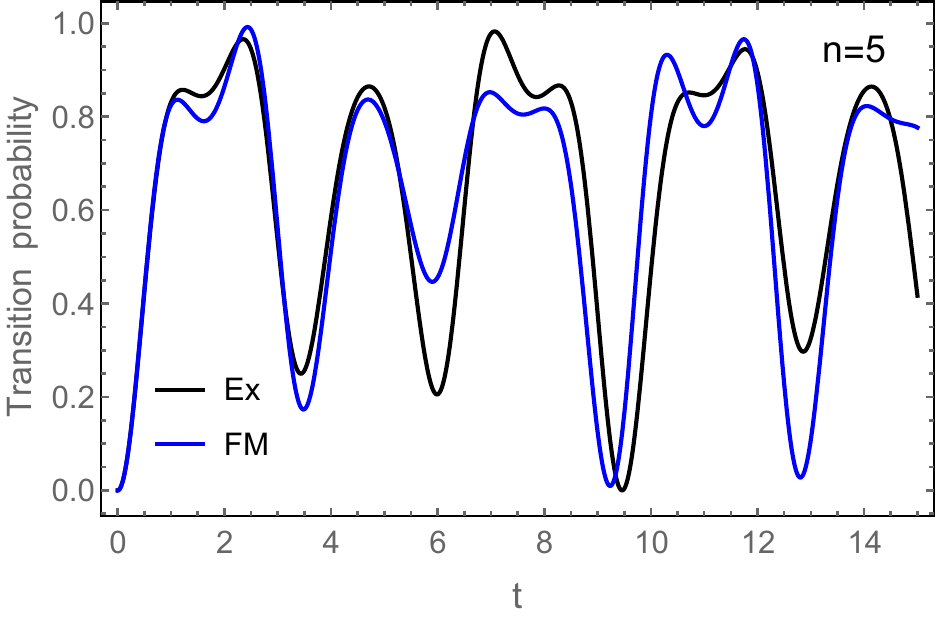}\\
	\includegraphics[scale=0.38]{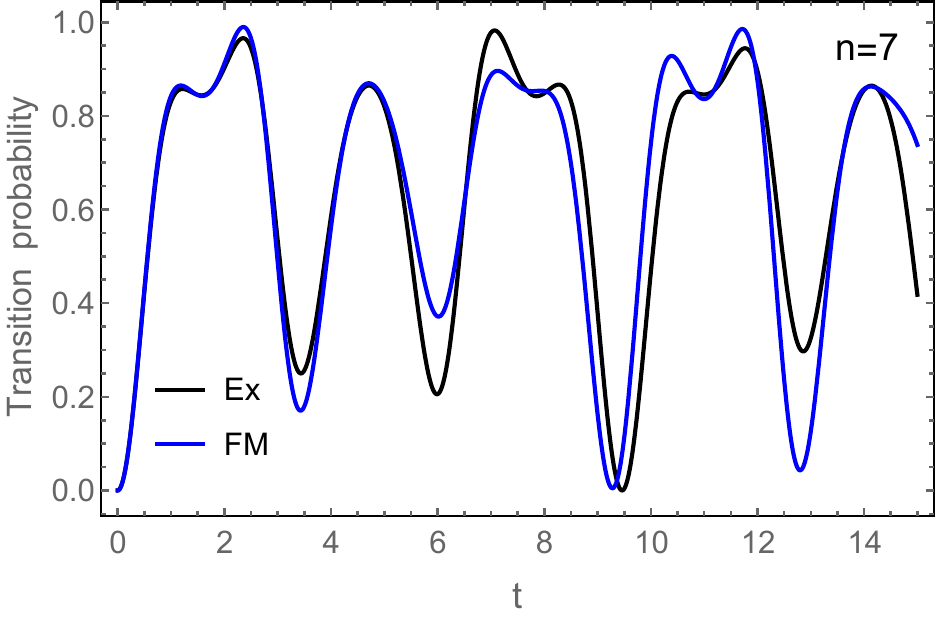}\includegraphics[scale=0.38]{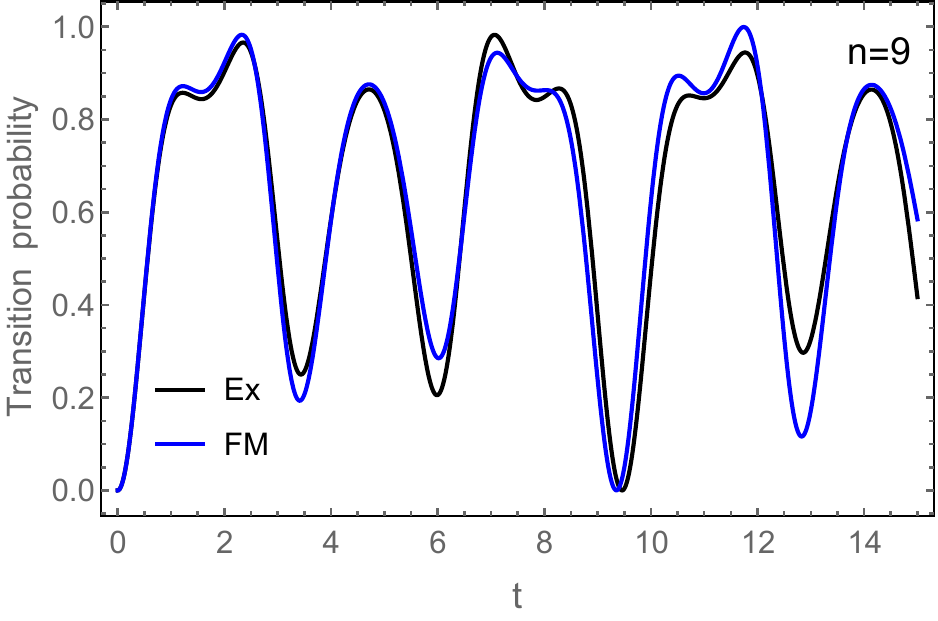}
	\caption{{\small Transition probability obtained with FM (solid blue line) with $n=3, 5, 7, 9$ terms (left to right, top to bottom)
			in comparison with the exact result (black line) for the Bloch--Siegert Hamiltonian with $\varepsilon = 1.5$.}}  \label{figu66}
\end{figure}

To analyze how incorporating an increasing number of terms in the Floquet--Magnus expansion leads to an improvement of the approximation, even for
large values of $\varepsilon$, we show in Figure \ref{figu66} the results for the transition probability for $\varepsilon = 1.5$ obtained with $n=3, 5, 7, 9$ terms
(left to right, top to bottom) as a function of time. Notice that even for this large value of the perturbation parameter, we still have a correct description of
the dynamics when a sufficiently large number of terms is taken into account.

\begin{figure}[H]
	\includegraphics[scale=0.49]{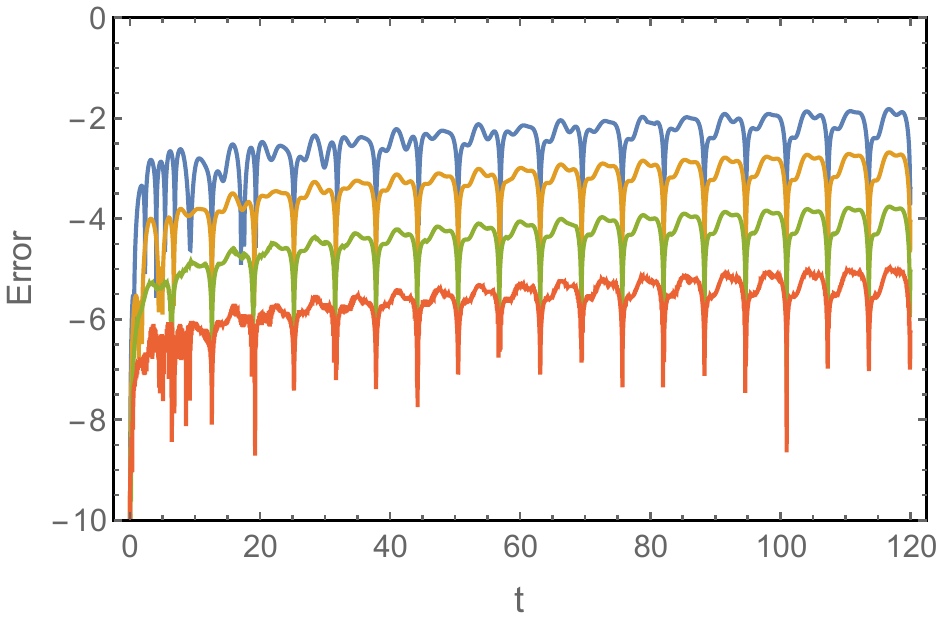} \includegraphics[scale=0.44]{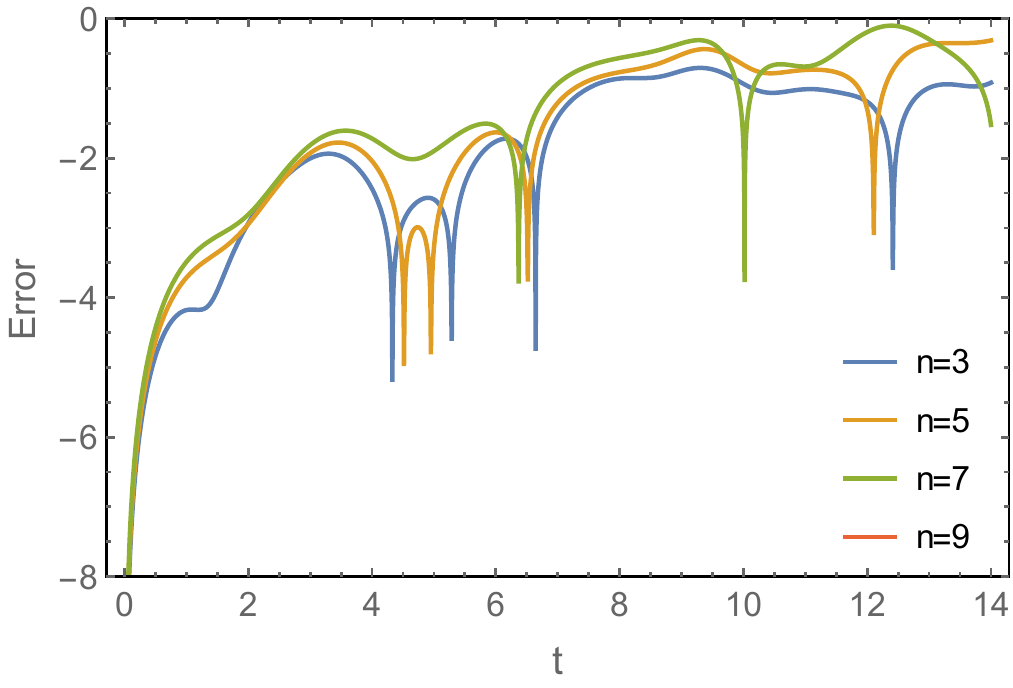}
	\caption{{\small Error (in logarithmic scale) in the transition probability committed by FM (left) and LD (right) with $n=3, 5, 7, 9$ terms each for the Bloch-Siegert
				Hamiltonian with $\varepsilon = 0.5$. Notice the different time interval in each panel.}} \label{figu77}
	\end{figure}

To get a more quantitative description of the different approximations, in Figure \ref{figu77} we depict the error in the transition probability 
committed by the Floquet--Magnus expansion (left) and Lie--Deprit (right) with $n=3, 5, 7, 9$ terms each for $\varepsilon = 0.5$. Notice that
taking into account more terms leads to a reduction in the error, and that the LD approximation clearly deteriorates with time, whereas FM still provides
excellent results even for large time intervals.

	%%%%%%%%%%%%%%%%%%%%%%%%%%%%%%%%%%%%%%%%%%

\section{Concluding remarks}

The subject of the perturbative treatment of linear systems of time-dependent differential equations has a long history in physics and mathematics.
In physics, in particular, it plays a fundamental role in the study of the evolution of quantum systems. Among the different procedures, the so-called
exponential perturbation algorithms have the relevant property of preserving the unitary character of the evolution. More in general, the approximations
they provide belong to the same Lie group than the exact solution of the problem. Archetypical examples of exponential perturbation theories
are the Magnus expansion since its inception in the 1950s, and more recently, the Floquet--Magnus expansion, several quantum averaging procedures
and a generalization of the well known Hori--Deprit perturbation theory of classical mechanics. Each of these algorithms has been derived 
in an independent way and it is not always easy to establish connections and common points between them.

The present paper tries to bridge this gap by showing that all of them can be seen in fact as the result of linear 
changes of coordinates expressed as the exponential of a certain 
series whose terms can be obtained recursively. In addition to the Magnus and Floquet--Magnus expansions, other techniques can also be
incorporated into our general framework, including the quantum averaging method and the Lie--Deprit perturbative algorithm. 
In the process we have also considered a novel approach, namely an exponential transformation rendering the original system into another
one in which the perturbation is removed. The resulting approximations preserve whatever qualitative features the exact solution may possess
(unitarity, orthogonality, symplecticness, etc.). 
Even the standard
perturbation theory in the interaction picture is recovered in this setting when the  exponential defining the transformation is truncated. 

With this same framework one might of course consider other exponential transformations, 
and this would automatically lead to new perturbation formalisms which might
be particularly adapted to the problem at hand. For instance, we could choose $F(t)$ as a diagonal time-dependent matrix, so that equation
(\ref{ccv.2}) is easy to solve, and construct the corresponding generator $\Omega(t)$.

We believe that the proposed framework sheds new light into nature of the different
expansions and allows one to create a unique procedure to analyze a given problem, obtaining all the expansions by choosing appropriately the coefficient
matrix in the new variables. Moreover, it also allows one to determine in an easy way what is the best algorithm for a given problem. This feature has been
illustrated in this work by applying the procedure to two simple examples. 

We have seen, in particular, that  for problems with a periodic time dependence, the Floquet--Magnus expansion leads to more accurate results
over longer time intervals, even for perturbed problems of the form $A(t) = A_0 + \varepsilon A_1(t)$ in the interaction picture,
unless the parameter $\varepsilon$ is exceedingly small,
in which other purely perturbative procedures are competitive. It would be interesting to check if this conclusion remains valid for more involved problems.

The interested reader may find useful the file available in \cite{arnal20wp} containing the \emph{Mathematica} implementation of the techniques exposed here for the quantum mechanical problem describing
the Bloch-Siegert dynamics.

In this paper we have always taken $X(0) = x(0)=x_0$, but it is clear that we can also take instead $X(0)$ as the image by the
transformation of $x_0$. In that case the solution
is factorized as 
\[
x(t) = \exp(\Omega(t)) \exp(t F) X(0) = \exp(\Omega(t)) \exp(t F) \exp(-\Omega(0)) x_0,
\]
where, obviously, the recurrences for determining $\Omega(t)$ and  $F$ are slightly different from those collected here. The final results can also 
vary, as shown  in \cite{arnal16apa}.

\vspace{1cm}

\appendix
\counterwithin*{equation}{section}
\renewcommand\theequation{\thesection\arabic{equation}}
\section{Appendix}

In this appendix we solve the equation
\begin{equation}   \label{ape.1}
\dot{\Omega}_n + [\Omega_n, A_0] = \mathcal{F}_n - F_n, 
\end{equation}
with initial condition $ \Omega_n(0) = 0$ and where $\mathcal{F}_n(t)$ is a quasi-periodic function
\[
\mathcal{F}_n(t) =   \sum_{k \in \mathbb{Z}^r} f_{n,k} \, \e^{i (k,  \omega) t},
\]
in such a way that $F_n$ is constant and $\Omega_n(t)$ is quasi-periodic with the same frequencies as $\mathcal{F}_n(t)$.  

First we compute the limiting mean value of $\mathcal{F}_n(t)$,
\[
\langle \mathcal{F}_n \rangle  \equiv \lim_{T\rightarrow\infty}\frac{1}{T} \int_0^{T} \mathcal{F}_n(t) dt
\]
and then the antiderivative
\begin{equation}   \label{emen}
M_n(t) = \int \e^{-t \ad_{A_0}} \left( \mathcal{F}_n(t) - \langle \mathcal{F}_n \rangle  \right) dt.
\end{equation}
Then, a direct computation show that 
\begin{equation}  \label{efen2}
F_n = \langle \mathcal{F}_n(t) \rangle - [A_0, M_n(0)], \qquad \Omega_n(t) = - M_n(0) + \e^{t \, \ad_{A_0}} \, M_n(t)
\end{equation}
provides the solution of (\ref{ape.1}). To show that $\Omega_n(t)$ is indeed quasi-periodic, we proceed as follows.
First, it is clear that
\begin{equation}  \label{efen3}
\mathcal{F}_n(t) - \langle \mathcal{F}_n \rangle = \sum_{k \in \mathbb{Z}^r \backslash \{0\}} f_{n,k} \, \e^{i (k,  \omega) t}
\end{equation}   
and thus
\[
\begin{aligned}
M_n(t) & =  \int dt \sum_{k \in \mathbb{Z}^r \backslash \{0\}}   \e^{-t \, \ad_{A_0}}  f_{n,k} \, \e^{i (k,  \omega) t}  = 
\int dt \sum_{k \in \mathbb{Z}^r \backslash \{0\}} \e^{-t \, \ad_{A_0}+i (k,\omega) t I}  f_{n,k} \\
& = \int dt \sum_{k \in \mathbb{Z}^r \backslash \{0\}} \e^{-t \, \ad_{A_0}+i (k,\omega) t I}  (-\ad_{A_0}+i (k,\omega) I ) g_{n,k}\,  \\
& = \sum_{k \in \mathbb{Z}^r \backslash \{0\}} \e^{-t \, \ad_{A_0}+i (k,\omega) t I}  g_{n,k} 
=   \e^{-t \, \ad_{A_0}} \sum_{k \in \mathbb{Z}^r \backslash \{0\}} g_{n,k} \, \e^{i (k,  \omega) t}, 
\end{aligned}
\]
where $g_{n,k}$ is the \textit{unique} solution of the linear system consisting of $d^2$ equations
\begin{equation}  \label{ape.2}
\Big(- \ad_{A_0}  + i (k, \omega) I \Big) X = f_{n,j}.
\end{equation}
In that case,
\begin{equation}  \label{omeapen}
\Omega_n(t) =  -M_n(0) + \sum_{k \in \mathbb{Z}^r \backslash \{0\}} g_{n,k} \, \e^{i (k,  \omega) t}
\end{equation}
is clearly a quasi-periodic function with the same basic frequencies as $\mathcal{F}_n(t)$.  On the other hand, the system (\ref{ape.2})
has a unique solution if and only if $X=0$ is the unique solution of
\[
\Big( -\ad_{A_0}  + i (k, \omega) I \Big) X  = -A_0 X +X (A_0 + i (k, \omega) I) = 0,
\]
and this happens when $A_0$ and $A_0 + i (k, \omega) I$
do not have common characteristic values, or equivalently when 
$\lambda_{\ell} - \lambda_m - i (k, \omega) \ne 0$ for all
$\ell, m \in \{1,\ldots, s\}$, where $\{\lambda_j\}_{j=1}^s$ denote the distinct eigenvalues of $A_0$. To avoid the problem of small divisors,
a diophantine
condition is introduced, namely one assumes that 
\begin{equation}  \label{diophantine}
|\lambda_{\ell} - \lambda_m - i (k, \omega)|  > \frac{\delta}{|k|^{\gamma}}    
\qquad \forall \; 
\ell, m \in \{1,\ldots, s\} \quad \forall \;  k \in \mathbb{Z}^r \backslash \{0\}
\end{equation}
for some constants $\delta>0$ and $\gamma > r-1$, see \cite{jorba97ero,braaksma87oaq,broer12raf} for more details. Here $|k| = |k_1| + \cdots +|k_r|$.

\vspace{2cm}

{\it This research has been funded by Ministerio de Econom\'{\i}a, Industria y 
Competitividad (Spain) through project MTM2016-77660-P (AEI/FE\-DER, UE)
and by Universitat Jaume I (projects UJI-B2019-17 and GACUJI/2020/05).}

\bibliographystyle{siam}

\end{document}